\documentclass{amsart}

\usepackage{amsmath}
\usepackage{amscd}
\usepackage{amsfonts}
\usepackage{amssymb}
\usepackage{latexsym}
\usepackage{pifont}
\usepackage{eepic}
\usepackage{array}

\newcommand{\ncm}{\newcommand}

\ncm{\beq}{\begin{equation}}
\ncm{\eeq}{\end{equation}}

\newtheorem{thm}{Theorem}[section]
\newtheorem{pro}[thm]{Proposition}
\newtheorem{lem}[thm]{Lemma}
\newtheorem{cor}[thm]{Corollary}
\newtheorem{thm&def}[thm]{Theorem \& Definition}

\theoremstyle{definition}
\newtheorem{defi}[thm]{Definition}

\theoremstyle{remark}
\newtheorem{rmk}[thm]{Remark}

\numberwithin{equation}{section}

\def\M{\mathsf{M}}
\ncm{\Cgd}{\mathsf{Cgd}}
\ncm{\BGD}{\mathsf{Bgd}}
\ncm{\BgdMap}{\mathsf{BgdMap}}
\ncm{\Bim}{\mathbb{BIM}}
\ncm{\YD}{\mathcal{YD}}
\ncm{\Gal}{\mathsf{Gal}}
\ncm{\Fib}{\mathsf{Fib}}
\ncm{\DHR}{\mathsf{DHR}}

\def\E{\mathcal{E}}
\ncm{\F}{\mathcal{F}}
\ncm{\G}{\mathcal{G}}

\ncm{\V}{\mathcal{V}}
\ncm{\W}{\mathcal{W}}
\ncm{\Z}{\mathcal{Z}}

\ncm{\asso}{\mathbf{a}}
\ncm{\luni}{\mathbf{l}}
\ncm{\runi}{\mathbf{r}}

\ncm{\End}{\operatorname{End}}
\ncm{\Hom}{\operatorname{Hom}}
\ncm{\BiEnd}{\operatorname{BiEnd}}

\def\Ker{\mbox{\rm Ker}\,}
\def\id{\mbox{\rm id}\,}

\newcommand{\ci}{\circ}

\def\o{\otimes}
\def\x{\times}

\ncm{\amalgo}[1]{\underset{\scriptscriptstyle #1}{\o}}
\ncm{\ractB}{\underset{\scriptscriptstyle B}{\ract}}
\ncm{\ractT}{\underset{\scriptscriptstyle T}{\ract}}
\ncm{\mash}{\Pisymbol{psy}{35}}
\ncm{\mashed}[1]{\underset{\scriptscriptstyle #1}{\Pisymbol{psy}{35}}}

\ncm{\oA}{\amalgo{A}}
\ncm{\oB}{\amalgo{B}}
\ncm{\oC}{\amalgo{C}}
\ncm{\oR}{\amalgo{R}}
\ncm{\oT}{\amalgo{T}}
\ncm{\oL}{\amalgo{L}}
\ncm{\oS}{\amalgo{S}}
\ncm{\oN}{\amalgo{N}}
\ncm{\oH}{\amalgo{H}}
\ncm{\oQ}{\amalgo{Q}}
\ncm{\ex}[1]{\underset{\scriptscriptstyle #1}{\x}}
\def\bra{\langle}
\def\ket{\rangle}
\ncm{\rarr}[1]{\stackrel{#1}{\longrightarrow}}
\ncm{\larr}[1]{\stackrel{#1}{\longleftarrow}}
\def\cop{\Delta}

\def\eps{\varepsilon}
\ncm{\op}{\mathrm{op}}
\ncm{\coop}{\mathrm{coop}}
\ncm{\co}{\mathrm{co}}
\ncm{\Fi}{\varphi}
\ncm{\OR}{\overrightarrow}
\ncm{\OL}{\overleftarrow}

\def\zeroT{^{(0)}}
\def\oneT{^{(1)}}
\def\twoT{^{(2)}}

\def\zeroB{_{(0)}}
\def\oneB{_{(1)}}
\def\twoB{_{(2)}}

\def\oneR{^{[1]}}
\def\twoR{^{[2]}}

\def\oneL{_{[1]}}
\def\twoL{_{[2]}}

\ncm{\coa}[1]{^{\langle #1\rangle}}
\ncm{\coabar}[1]{^{\langle \overline{#1}\rangle}}

\ncm{\I}{\mathcal{I}}

\def\du1{\hat 1}
\def\iso{\stackrel{\sim}{\rightarrow}}

\def\ract{\triangleleft}

\ncm{\under}{\mbox{\rm\_}\,}
\ncm{\Cnt}{\mathsf{C}}

\ncm{\ZZ}{\mathbb{Z}}

\ncm{\Oo}{\mathcal{O}}
\ncm{\Ha}{\mathcal{H}}
\ncm{\Ve}{\mathcal{V}}

\ncm{\sqrM}{$\sqrt{\text{Morita}}$}
\ncm{\ld}{\,.\,}
\ncm{\ud}{\,^.\,}
\ncm{\Mor}{\underset{k}{\sim}}
\ncm{\into}{\hookrightarrow}
\ncm{\coinv}[1]{{\co\text{-}#1}}
\ncm{\lZ}{\overrightarrow{\mathcal{Z}}}
\ncm{\rZ}{\overleftarrow{\mathcal{Z}}}
\ncm{\fgp}{\mathrm{fgp}}
\ncm{\adj}{\dashv}
\ncm{\Mref}{{\scriptscriptstyle M\text{-}\mathrm{ref}}}

\begin{document}

\title{Finitary Galois extensions over noncommutative bases}
\author[I. B\'alint, K. Szlach\'anyi]{Imre B\'alint \and Korn\'el Szlach\'anyi}
\address{Research Institute for Particle and Nuclear Physics, Budapest}
\email{balint@rmki.kfki.hu, szlach@rmki.kfki.hu}
\thanks{Supported by the Hungarian Scientific Research Fund, OTKA T-034 512
and T-043 159}

\begin{abstract}
We study Galois extensions $M^{(\co\text{-})H}\subset M$ for 
$H$-(co)module algebras $M$ if $H$ is a Frobenius Hopf algebroid.
The relation between the action and coaction pictures is analogous to that found
in Hopf-Galois theory for finite dimensional Hopf algebras over fields. 
So we obtain generalizations of various classical theorems of 
Kreimer-Take\-uchi, Doi-Takeuchi and Cohen-Fischman-Montgomery. 
We find that the Galois extensions $N\subset M$ over some Frobenius Hopf 
algebroid are precisely the balanced depth 2 Frobenius extensions. 
We prove that the Yetter-Drinfeld categories over $H$ are always braided and 
their braided commutative algebras play the role of noncommutative scalar 
extensions by a slightly generalized Brzezi\'nski-Militaru Theorem. 
Contravariant "fiber functors" are used to prove an analogue of Ulbrich's 
Theorem and to get a monoidal embedding of the module category $\M_E$ of the 
endomorphism Hopf algebroid $E=\End\,_NM_N$ into $_N\M_N^\op$.
\end{abstract}

\maketitle

\section{Introduction}

The problem of extending Hopf Galois theory to quantum groupoids has been 
attracting some attention in recent years. That this theory should possess  
interesting new applications even for finite quantum groupoids is manifest 
already from the pioneering work of D. Nikshych and L. Vainerman \cite{N-V}. 
A pure algebraic Galois theory for weak Hopf algebras has been proposed 
by S. Caenepeel and E. de Groot \cite{Caenepeel-DeGroot}. 
As the next step of generalization, this paper is devoted to 
developing a Galois theory for Frobenius Hopf algebroids. These quantum 
groupoids are the analogues of finite dimensional Hopf algebras over a 
field or Frobenius Hopf algebras over a commutative ring. Therefore it is not 
surprising that we obtain generalizations of the classical theorems of 
Kreimer-Takeuchi \cite{Kreimer-Takeuchi}, Doi-Takeuchi \cite{Doi-Takeuchi} and 
Cohen-Fischman-Montgomery \cite{CFM} (see Theorems \ref{thm: def Gal} and 
\ref{thm: WS}). Our results partly overlap with those of the recent paper 
\cite{B: Gal} by G. B\"ohm who studies Galois theory for general Hopf 
algebroids using previous results from the theory of corings 
\cite{Brz: corings,Brz-Wis}. In our approach the double algebraic structure 
\cite{Sz: DA} of Frobenius Hopf algebroids is particularly useful e. g. in 
proving that Yetter-Drinfeld categories are braided (Proposition \ref{pro: YD 
braided}), in this way generalizing a result of \cite{Caenepeel-Wang-Yin}, or in 
obtaining an instrinsic characterization of Galois extensions as being the depth 
2, balanced, Frobenius extensions (Theorem \ref{thm: bD2F}).

\subsection{Modules and comodules over bialgebroids}

Let $k$ be a commutative ring. We choose  the category $\M=\M_k$ of $k$-modules 
as our base category. This means that all objects and morphisms we use
have an underlying $k$-module or $k$-module morphism, respectively.
In particular, algebras are always meant to be $k$-algebras and unadorned $\o$ 
means tensor product in $\M_k$.

Let $T$ be an algebra and let $T^e:=T^\op\o T$ be its enveloping algebra. 
A right bialgebroid over $T$ consists of
\begin{itemize}
\item an algebra $A$
\item a $T^e$ ring structure on $A$, i.e., an algebra morphism 
$t_r\o s_r:T^e\to A$
\item and a $T$-coring structure $\bra A_{T^e},\cop_T,\Fi_T\ket$
\end{itemize}
subject to axioms, see e.g. \cite{K-Sz}. 

If $A$ is a right bialgebroid over $T$ then
a right module over $A$ is the same thing as a right module over the $k$-algebra 
$A$ and a right $A$-module map is defined accordingly. The $T^e$-ring structure 
$T^e\to A$ endows the category $\M_A$ of right $A$-modules with a (monadic) 
forgetful functor $U:\M_A\to\,_T\M_T$ by identifying $\M_{T^e}$ with $_T\M_T$.
The coring structure of $A$ serves to make $\M_A$ a monoidal category. The 
monoidal product of the $A$-modules $V$ and $W$ is the $k$-module $V\oT W$ 
together with the right $A$-action $(v\oT w)\ract a:=(v\ract a\oneT)\oT(w\ract 
a\twoT)$. In this way the forgetful functor $U$ becomes strict monoidal. 

Left bialgebroids and their category of left modules can be defined by passing 
to the opposite algebra in all occurences of an algebra in the definition of a 
right bialgebroid and their right modules. So let $B$ be an algebra which stands 
for $T^\op$ and let $B^e:=B\o B^\op$. Then a left bialgebroid over $B$ consists 
of
\begin{itemize}
\item an algebra $A$
\item a $B^e$-ring structure on $A$, i.e., an algebra morphism 
$s_l\o t_l:B^e\to A$
\item and a $B$-coring structure $\bra\,_{B^e}A,\cop_B,\Fi_B\ket$.
\end{itemize}
The category of left $A$-modules has a monoidal product $V\oB W$ such that the 
forgetful functor $_AV\mapsto \,_{B^e}V\equiv\,_BV_B$ is strict monoidal.

Right comodules can be defined for both left and right bialgebroids as follows.
Let $A$ be a right bialgebroid over $T$. Then a right $A$-comodule consists of
\begin{itemize}
\item a right $T$-module $X$
\item a right $T$-module map $\delta:X\to X\oT A$
\end{itemize}
such that
\begin{align*}
(\delta\oT A)\circ\delta&=(X\oT\cop_T)\circ\delta\\
(X\oT\Fi_T)&=X
\end{align*}
suppressing the coherence isomorphisms of $_T\M_T$.
A morphism of comodules $\tau:\bra X,\gamma\ket\to \bra Y,\delta\ket$ is a right 
$T$-module map $\tau:X\to Y$ satisfying $(\tau\oT 
A)\circ\gamma=\delta\circ\tau$. The category of right $A$-comodules is denoted 
$\M^A$.

The above definition of comodules disguises the fact that $\M^A$ is 
monoidal with a strict monoidal forgetful functor $\M^A\to\,_T\M_T$.
Notice that although $M$ is not a left $T$-module, $M\oT A$ is by setting
$t\cdot(x\oT a)=x\oT s_r(t)a$.
\begin{pro} \label{pro: coactions are bim}
Let $\bra X,\delta\ket$ be a right comodule over the rigt bialgebroid $A$.
Then $X$ has a unique left $T$-module structure such that $\delta$ is a left 
$T$-module map. With this left module structure
\begin{enumerate}
\item $X$ is a $T$-$T$-bimodule,
\item $\delta$ is a $T$-$T$-bimodule map,
\item $\delta(X)\subset X\x_T A$,
\item and every arrow $\tau\in\M^A$ is a $T$-$T$-bimodule map.
\end{enumerate}
In (3) we used Takeuchi's $\x$-product which is defined by
\begin{equation*}
X\x_T A:=\{\sum_i x_i\oT a_i\in X\oT A\,|\,\sum_i t\cdot 
x_i\oT a_i=\sum_i x_i\oT t_r(t)a_i\,\forall t\in T\}.
\end{equation*}
\end{pro}
\begin{proof}
If $X$ is a left $T$-module and $\delta$ is a left $T$-module map then
\begin{align*}
t\cdot x&=(t\cdot x)\zeroT\cdot\Fi_T((t\cdot x)\oneT)\\
&=x\zeroT\cdot\Fi_T(s_r(t)x\oneT)
\end{align*}
This proves uniqueness. If we use the above formula to define $t\cdot x$ then we 
find that it is a left action because $s_r:T\to A$ is an algebra 
homomorphism. It commutes with the right $T$-action 
\begin{equation*}
t\cdot (x\cdot t')=x\zeroT\cdot\Fi_T(s_r(t)as_r(t'))=(t\cdot x)\cdot t'
\end{equation*}
so $X$ is a $T$-$T$-bimodule and the coaction is a bimodule map,
\begin{align*}
\delta(t\cdot x\cdot t')=x\zeroT\oT s_r(t)x\oneT s_r(t').
\end{align*}
Now the Takeuchi property (3) holds automatically,
\begin{align*}
t\cdot x\zeroT\oT x\oneT&=x\zeroT\cdot\Fi_T(s_r(t)x\oneT)\oT x\twoT\\
&=x\zeroT\oT\Fi_T(t_r(t)x\oneT)\cdot x\twoT\\
&=x\zeroT\oT t_r(t)x\oneT\,.
\end{align*}
If $\tau:X\to Y$ is a comodule morphism then
\begin{align*}
\tau(t\cdot x)&=\tau(x\zeroT)\cdot\Fi_T(s_r(t)x\oneT)=
\tau(x)\zeroT\cdot\Fi_T(s_r(t)\tau(x)\oneT)\\
&=t\cdot\tau(x)\,.
\end{align*}
\end{proof}

The tensor product of right comodules $X$ and $Y$ can now be defined as $X\oT Y$ 
with coaction
\begin{equation}\label{tensor prod of comodules over right bgd}
(x\oT y)\zeroT\oT (x\oT y)\oneT=(x\zeroT\oT y\zeroT)\oT x\oneT y\oneT\,.
\end{equation}
This makes the category of right $A$-comodules $\M^A$ monoidal and the 
forgetful functor $\M^A\to\,_T\M_T$ strict monoidal.

For left bialgebroids $A$ over $B$ a right comodule is an arrow $\delta_A:M\to 
M\oB A\in\M_B$ satisfying coassociativity and counitality. A right comodule 
carries a left $B$-module structure such that $\delta_A$ is a $B$-$B$-bimodule 
map and such that $\M^A$ is a monoidal category with strict monoidal forgetful 
functor to $_B\M_B$. The monoidal product of two right comodules $X$ and $Y$ 
has coaction
\begin{equation}\label{tensor prod of comodules over left bgd}
(x\oB y)\zeroB\oB (x\oB y)\oneB=(x\zeroB\oB y\zeroB)\oB y\oneB x\oneB\,.
\end{equation}
Note the different order compared to (\ref{tensor prod of comodules over 
right bgd}).

\subsection{Double algebras}

Studying module (co)algebras over bialgebroids one can obtain  
generalizations of certain theorems of Hopf-Galois theory and, 
therefore, hints toward the proper definition of bialgebroid Galois extensions
\cite{K-Sz,Sz: Strasbourg,Kadison}. For example, the behaviour of depth 2 
balanced extensions $N\subset M$ of algebras suggest that they are 
Galois extensions in the very noncommutative sense. Of course, in the absence of 
antipode, even in the finitely generated projective (fgp) case, many results of 
classical and Hopf Galois theory are far from reach.

This leads us to study Hopf algebroids instead and for finiteness we assume that 
it has an integral which is a Frobenius functional. These Frobenius Hopf 
algebroids were shown in \cite{Sz: DA} to possess a \textit{distributive double 
algebra} structure (DDA) by choosing a Frobenius integral. Here we summarize 
its basic properties. 

A double algebra is a $k$-module $A$ equipped with two associative unital 
multiplications: the vertical multiplication, denoted $a\ci a'$, has unit 
element $e$ and the horizontal multiplication, denoted $a\star a'$, has unit 
element $i$. So we have the horizontal and vertical algebras $H=\bra 
A,\star,i\ket$ and $V=\bra A,\ci,e\ket$, respectively. The multiplications with 
the wrong unit, i.e.,
\begin{alignat*}{2}
\Fi_L(a)&:=a\star e &\qquad\Fi_R(a)&:=e\star a \\
\Fi_B(a)&:=a\ci i &\qquad \Fi_T(a)&:=i\ci a 
\end{alignat*}
map onto subalgebras $L$ and $R$ of $V$ and $B$ and $T$ of $H$.
Assuming for $X=L,R,B,T$ that the algebra extensions $X\subset A$ are Frobenius 
with Frobenius homomorphism $\Fi_X$ we obtain the notion of Frobenius DA's.
In this way $A$ has Frobenius algebra structures in all the bimodule 
categories $_X\M_X$ for $X=L,R,B,T$ which implies four 
comultiplications

$\bra A, \Delta_B,\Fi_B\ket$ is a comonoid in $_B\M_B$, where
      $\Delta_B(a)\equiv a\oneB\oB a\twoB=a\star u_k\oB v_k$,

$\bra A, \Delta_L,\Fi_L\ket$ is a comonoid in $_L\M_L$, where
      $\Delta_L(a)\equiv a\oneL\oL a\twoL=a\ci x_j\oL y_j$,

$\bra A, \Delta_T,\Fi_T\ket$ is a comonoid in $_T\M_T$, where
      $\Delta_T(a)\equiv a\oneT\oT a\twoT=a\star u^k\oT v^k$,

$\bra A, \Delta_R,\Fi_R\ket$ is a comonoid in $_R\M_R$, where
      $\Delta_R(a)\equiv a\oneR\oR a\twoR=a\ci x^j\oR y^j$.

\noindent
where note the special notation for the dual bases of the base homomorphisms 
$\Fi_X$. It turns out \cite[Proposition 3.2]{Sz: DA} that vertical 
multiplication with the horizontal type of comultiplications $\cop_B$ and 
$\cop_T$ obey bialgebroid like relations. However, if we also postulate the 
distributivity rules
\begin{align}
a\ci(a'\star a'')&=(a\oneB\ci a')\star(a\twoB\ci a'') \label{eq: DB}\\
a\star(a'\ci a'')&=(a\oneL\star a')\ci(a\twoL\star a'') \label{eq: DL}\\
(a'\star a'')\ci a&=(a'\ci a\oneT)\star(a''\ci a\twoT) \label{eq: DT}\\
(a'\ci a'')\star a&=(a'\star a\oneR)\ci(a''\star a\twoR) \label{eq: DR}
\end{align}
in which case we say that $\bra A,\ci,e,\star,i\ket$ is a distributive double 
algebra (DDA), then $V$ and $H$ become Hopf algebroids \cite{B-Sz: Hgd} in 
duality. The underlying left bialgebroids are
\[
\bra V,B,\Fi_L|_B,\Fi_R|_B,\Delta_B,\Fi_B\ket \quad\text{and}\quad
\bra H,L,\Fi_B|_L,\Fi_T|_L,\Delta_L,\Fi_L\ket
\]
and the right bialgebroids are 
\[
\bra V,T,\Fi_R|_T,\Fi_L|_T,\Delta_T,\Fi_T\ket \quad\text{and}\quad
\bra H,R,\Fi_T|_R,\Fi_B|_R,\Delta_R,\Fi_R\ket
\]
The notation means e.g. that $V$ over $T$ has source map $s_r:t\mapsto 
\Fi_R(t)$, target map $t_r:t\mapsto \Fi_L(t)$ and counit $\Fi_T$. Or, $H$ over 
$R$ has source map $s_r:r\mapsto \Fi_T(r)$, target map $t_r:r\mapsto \Fi_B(r)$, 
and counit $\Fi_R$. The antipode of $V$ -- called the antipode of the double 
algebra -- is an antiautomorphism $S$ which is also an antiautomorphism of $H$ 
but the antipode of $H$ is $S^{-1}$. (There is a regrettable mistake in 
\cite[Theorem 7.4]{Sz: DA} where $H$ was claimed to have antipode also $S$; see 
\texttt{arXiv: math.QA/0402151 v2} for the corrected version.) The vertical Hopf 
algebroid has Frobenius integral $i$ and $H$ has $e$.

\section{Modules and comodules over DDA's}

\subsection{Modules}
Let $\bra A,\ci,e,\star,i\ket$ be a double algebra. A \textit{right $A$-module} 
is a $k$-module $M$ together with an associative unital action $M\oR H\to M$ of 
the horizontal algebra $H=\bra A,\star,i\ket$ denoted $m\oR h\mapsto m\ract 
h$.

Equivalently, a right $A$-module can be formulated in the category $\M_{B\o 
T}$ as an object $M_{B\o T}$ and an arrow $M\amalgo{B\o T}A\to M$ satisfying 
associativity and unitality w.r.t the algebra $H$ in $_{B\o T}\M_{B\o T}$. The 
$T$ and $B$-actions are denoted by $m\ud t$ and $m\ld b$, respectively.

Analogously one can define \textit{left $A$-modules} as left modules over $H$ 
and \textit{bottom} and \textit{top $A$-modules} as "left", respectively 
"right", modules over the vertical algebra $V=\bra A,\ci,e\ket$.

\subsection{Comodules}

A \textit{right $A$-comodule} over a Frobenius double algebra consists of an 
object $M$ and two arrows $\delta_M\colon M\to M\oB A$, $\delta^M\colon M\to 
M\oT A$ in $\M_{B\o T}$ such that
\begin{itemize}
\item $\bra M_B,\delta_M\ket$ is a right comodule over the left bialgebroid $V$ 
over $B$, 
\item $\bra M_T,\delta^M\ket$ is a right comodule over the right 
bialgebroid $V$ over $T$
\item and the two coactions satisfy the mixed coassociativity 
conditions
\begin{align} \label{mixed coa 1}
{m\zeroT}\zeroB\oB {m\zeroT}\oneB\oT m\oneT
&=m\zeroB\oB {m\oneB}\oneT\oT {m\oneB}\twoT\\
{m\zeroB}\zeroT\oT {m\zeroB}\oneT\oB m\oneB
&=m\zeroT\oT {m\oneT}\oneB\oB {m\oneT}\twoB \label{mixed coa 2}
\end{align}
\end{itemize}
where we used the notation
\begin{align*}
\delta_M(m)&=m\zeroB\oB m\oneB\\
\delta^M(m)&=m\zeroT\oT m\oneT
\end{align*}
for $m\in M$.

A \textit{right $A$-comodule morphism} $\tau:X\to Y$ is a right  $B\o T$-module 
map  
which is a right comodule morphism for both the left bialgebroid $V_B$ and the 
right bialgebroid $V_T$. The category of right $A$-comodules is denoted by 
$\M^V$.
The occurence of two compatible coactions in the definition of an $A$-comodule 
is precisely what we need to identify $\M^V$ and $\M_H$ in case of DDA's. 

\begin{lem} \label{lem: act-coact}
Let $A$ be a DDA and let
$\delta_M$ and $\delta^M$ be two coactions of $V_B$, respectively $V_T$, on 
$M$. They then determine two right $H$-actions on $M$,
\begin{align}
m\ractB h&=m\zeroB \ld\Fi_B(m\oneB\star h) \label{eq: coa->act B}\\
m\ractT h&=m\zeroT \ud\Fi_T(m\oneT\star h)\,. \label{eq: coa->act T}
\end{align}
The two actions coincide if and only if the two coactions satisfy the mixed 
coassociativity condition (\ref{mixed coa 1}) and (\ref{mixed coa 2}).
\end{lem}
\begin{proof}
The inverses of (\ref{eq: coa->act B}) and (\ref{eq: coa->act T}) can be given 
in terms of the dual bases of $\Fi_B$ and $\Fi_T$ as
\begin{align}
m\zeroT\oT m\oneT&=m\ractT u^k\oT v^k \label{eq: act->coa B}\\
m\zeroB\oB m\oneB&=m\ractB u_k\oT v_k \label{eq: act->coa T}
\end{align}
Therefore if $\ractB=\ractT$ then
\begin{align*}
{m\zeroT}\zeroB\oB {m\zeroT}\oneB\oT m\oneT
&=(m\ractT u^k)\ractB u_l\oB v_l\oT v^k
=m\ractB(u^k\star u_l)\oB v_l\oT v^k\\
&=m\ractB u_l\oB v_l\star u^k\oT v^k
=m\zeroB\oB {m\oneB}\oneT\oT {m\oneB}\twoT
\end{align*}
and similarly for (\ref{mixed coa 2}). On the other hand, if mixed 
coassociativity holds then
\begin{align*}
m\ractT h&=(m\ractT h)\zeroB\ld\Fi_B((m\ractT h)\oneB)
={m\zeroT}\zeroB\ld\Fi_B({m\zeroT}\oneB\star\Fi_T(m\oneT\star h))\\
&=m\zeroB\ld\Fi_B({m\oneB}\oneT\star\Fi_T({m\oneB}\twoT\star h))
=m\zeroB\ld\Fi_B(m\oneB\star h)\\
&=m\ractB h
\end{align*}
\end{proof}

If $M$ is a right module over the DDA $A$ then it is a right $V$-comodule $M^V$ 
and a right $H$-module $M_H$ at the same time. The invariants of $M_H$, 
\begin{align}
M^H&{:=}\{ n\in M | n\ract h=n\ract\Fi_T\Fi_R(h),\ h\in H\}\\
&=\{ n\in M | n\ract h=n\ract\Fi_B\Fi_R(h),\ h\in H\}
\end{align}
and the coinvariants of $M^V$,
\begin{align}
M^\coinv{V}&{:=}\{n\in M|n\zeroT\oT n\oneT=n\oT e\}\\
&=\{n\in M|n\zeroB\oB n\oneB=n\oB e\},\notag
\end{align}
yield one and the same $k$-submodule of $M$. This is an instance of the more 
general identification between the categories of $H$-modules, $V_B$-comodules,
and $V_T$-comodules. Since $\Fi_T$ and $\Fi_B$ restrict to 
algebra isomorphisms $R\to T$ and $R^\op\to B$, respectively \cite[Lemma 
2.2]{Sz: DA}, the identifications between $H$-modules and $V$-comodules provide 
a monoidal category isomorphism $\M^{V_T}\cong\M_H$ and the antimonoidal 
category isomorphism $\M^{V_B}\cong\M_H$.  We can use these isomorphisms to 
introduce $\oR$ both in $\M^{V_T}$ and $\M^{V_B}$ as the monoidal product while 
keeping $\oT$ and $\oB$ to appear in the coactions. One  advantage of this 
convention is that the difference between (\ref{tensor prod of comodules over 
left bgd}) and (\ref{tensor prod of comodules over right bgd}) disappears, viz. 
(\ref{eq: comalg T}) and (\ref{eq: comalg B}). Now the $R$ becomes a monoidal 
unit in three senses: As a right ideal in $H$ it is the trivial right 
$H$-module, $r\ract h=r\star h$. But it is also a right comodule over $V_T$ via 
$r\zeroT \oT r\oneT=e\oT r$ and a right comodule over $V_B$ via $r\zeroB \oB 
r\oneB=e\oB r$.

\subsection{Module (co)algebras} \label{ss: ext}

Comodule algebras over $V$ are monoids in $\M^V$ and therefore they are the same 
as monoids in $\M_H$, i.e., module algebras over $H$.

Hence a right $H$-module algebra $M$ consists of an algebra map 
$\eta:R\to M$ inducing the bimodule structure $_RM_R$ and a bimodule map 
$\mu:M\oR M\to M$, $m\oR m'\mapsto mm'$, satisfying $(mm')\ract h=(m\ract 
h\oneR)(m\ract h\twoR)$.
In the language of the $V$-coactions (\ref{eq: act->coa B}), 
(\ref{eq: act->coa T}) these correspond to the right comodule algebra relations
\begin{alignat}{2}
(mm')\zeroT\oT(mm')\oneT&=
m\zeroT{m'}\zeroT\oT m\oneT\ci{m'}\oneT
&\quad
1\zeroT\oT 1\oneT&=1\oT e   \label{eq: comalg T}\\
(mm')\zeroB\oB(mm')\oneB&=
m\zeroB{m'}\zeroB\oB m\oneB\ci{m'}\oneB
&\quad
1\zeroB\oB 1\oneB&=1\oB e   \label{eq: comalg B}
\end{alignat}
respectively. 
Just as in the case of Hopf algebras the invariants of a module algebra form a 
subalgebra. More precisely we have the following
\begin{lem} \label{lem: R->M}
For any right $H$-module $M$ there is a unique $k$-module map \newline
$\Hom_H(R,M)\to M^H$ that makes the diagram
\[
\begin{CD}
\Hom_H(R,M)@>\Hom(\Fi_R,M)>>\Hom_H(H,M)\\
@VVV @VV{f\mapsto f(i)}V\\
M^H@>\subset>>M
\end{CD}
\]
commutative. This $k$-module map is an isomorphism. If $M_H$ is a module algebra 
then the diagram is in the category of $k$-algebras. In particular, $M^H\subset 
M$ is a subalgebra which is isomorphic to the convolution algebra $\Hom_H(R,M)$. 
\end{lem}

The smash product $H\mash M$ for a right $H$-module algebra is defined to be
to $k$-module $H\oR M$ equipped with multiplication
\begin{equation}
(h\mash m)(h'\mash m')=h\star {h'}\oneR\mash (m\ract {h'}\twoR)m'
\end{equation}
and unit element $i\mash 1$.

Next we consider extensions. Let $N\to M^H\subset M$ be an algebra map. Then we 
have left actions $\lambda$ of $N$ and $\lambda$ of $M^H$ on $M$. Denoting
$\E:=\End(\,_NM)$ we have an algebra map $H\mash M\to\E$ by 
\begin{equation} \label{eq: Gamma}
m'\cdot(h\mash m):=(m'\ract h)m
\end{equation}
so that $M$ becomes an $N$-$(H\mash M)$-bimodule. We have the inclusions
\begin{equation} \label{eq: to be bal}
\lambda(N)\subset \End_\E(M)\subset \End_{H\mash M}(M)=\lambda(M^H)
\end{equation}
where the last equality can be proven exactly as in the Hopf algebra case 
\cite[8.3.2]{Montgomery}.
\begin{defi}
An algebra homomorphism $\eta:N\to M$ is called a right $A$-extension for some 
DDA $A$ if $M$ is a right module algebra over $A$ and $\eta$ factorizes through 
$M^H\subset M$ via an algebra isomorphism $N\iso M^H$.
\end{defi}
Later on an $A$-extension will be meant in the narrower sense that $N=M^H$ but 
sometimes, as in Section \ref{s: fiber} we need this more categorical 
definition.
\begin{lem} \label{lem: ext bal}
Let $A$ be a DDA and $N\to M$ be a right $A$-extension. Then 
\begin{enumerate}
\item $_NM$ is balanced, i.e., $\BiEnd(\,_NM)=\lambda(N)$ and
\item $_NM_{H\mash M}$ is faithfully balanced iff the canonical map $H\mash 
M\to\E$ given by (\ref{eq: Gamma}) is an isomorphism. 
\end{enumerate}
\end{lem}
\begin{proof}
Both statements are immediate consequences of the fact that all the inclusions 
in (\ref{eq: to be bal}) reduce to equalities in case of $A$-extensions.
\end{proof}

\section{Galois extensions}

\subsection{The coaction picture}

Let $M$ be a right comodule algebra over the Hopf algebroid $V$ and let 
$N:=M^\coinv{V}$. Then the maps
\begin{align}
\gamma^M&:M\oN M\to M\oT V\,,\qquad m\oN m'\mapsto m{m'}\zeroT\oT {m'}\oneT\\
\gamma_M&:M\oN M\to M\oB V\,,\qquad m\oN m'\mapsto m\zeroB m'\oB m\oneB
\end{align}
are $M$-$M$-bimodule maps if we endow $M\oT V$ and $M\oB V$ with the structure
\begin{align}
m'\cdot(m\oT v)\cdot m''&=m'm{m''}\zeroT\oT v\ci {m''}\oneT\,,\\
m'\cdot(m\oB v)\cdot m''&=m'\zeroB mm''\oB m'\oneB\ci v\,,
\end{align}
respectively. They are also right $V$-comodule maps, i.e., belong to $\M^V$, 
because they can be written as composites of $\mu_M$ and $\delta^M$, 
respectively $\mu_M$ and $\delta_M$.

\begin{lem} \label{lem: 2 gamma}
Let $M$ be a right $V$-comodule algebra over the Hopf algebroid $V$. Then
$\gamma^M$ is epimorphism iff $\gamma_M$ is and $\gamma^M$ is 
isomorphism iff $\gamma_M$ is.
\end{lem}
\begin{proof}
Let $\phi$ denote the composite 
\begin{align}
&\begin{CD}
M\oT V@>\delta_M\oT V>>M\oB V\oT V@>M\oB V\o S>>M\oB V\oR V
@>M\oB\mu_V>>M\oB V
\end{CD}\\
&\ \ \ \ m\oT v\mapsto m\zeroB\oB m\oneB S(v) \notag
\end{align}
where $S$ is the antipode of the Hopf algebroid $V$. Then $\phi$ has inverse
\[
\phi^{-1}(m\oB v)=m\zeroT\oT S^{-1}(v) m\oneT\,.
\]
and one obtains that $\phi\ci\gamma^M=\gamma_M$.
\end{proof}

The next result is an immediate generalization of \cite[Theorem 
8.3.1]{Montgomery}.
\begin{pro} \label{pro: epi implies iso}
Assume that $V$ is a Frobenius Hopf algebroid and $M$ is a right $V$-comodule 
algebra with coinvariant subalgebra $N$. Then $\gamma^M$ being epi implies that 
$\gamma^M$ is an isomorphism and $M_N$ is finitely generated projective.
\end{pro}
\begin{proof}
Let $V$ and $H$ be the vertical and horizontal Hopf algebroid of a distributive 
double algebra $\bra A,\ci,e,\star,i\ket$. Then $M$ is a right $H$-module 
algebra and $e$, the unit of $V$, is an integral for $H$, therefore $m\ract e\in 
N$, $m\in M$. By the hypothesis there exists $\sum_j m_j\oN m'_j\in M\oN M$ such 
that 
\[
\sum_j m_j{m'}_j\zeroT\oT {m'}_j\oneT=1\oT i\,.
\]
Therefore we can write for arbitrary $m\in M$ that
\begin{align*}
\sum_j m_j((m'_jm)\ract e)&=\sum_j m_j(m'_j\ract e\oneR)(m\ract e\twoR)\\
&=\sum_j m_j\left({m'}_j\zeroT\ud\Fi_T({m'}_j\oneT\star e\oneR)\right)(m\ract 
e\twoR)\\
&=(1\ract (i\ci e\oneR))(m\ract e\twoR)=(1\ract i\oneR)(m\ract i\twoR)\\
&=m\ract i\ =\ m
\end{align*}
proving that $(m'_j\under)\ract e$ is a dual basis of $m_j$ for $M_N$, thus 
$M_N$ is fgp.

Next we show that $\gamma_M$ is mono. Suppose $\sum_i z_i\oN 
w_i\in\Ker\gamma_M$. Then 
\[
\sum_i {z_i}\zeroB w_i\oB{z_i}\oneB\ =\ 0\,.
\]
Using the dual bases for $M_N$ we find that
\begin{align*}
\sum_i z_i\oN w_i&=\sum_i\sum_j m_j((m'_jz_i)\ract e)\oN w_i\\
&=\sum_j m_j\oN\sum_i\left({m'_j}\zeroB{z_i}\zeroB\ld\Fi_B(({m'_j}\oneB\ci 
{z_i}\oneB)\star e)\right)w_i\\
&=\sum_j m_j\oN\sum_i {m'_j}\zeroB{z_i}\zeroB w_i\ld\Fi_B({m'_j}\oneB\ci 
{z_i}\oneB)\\
&=0\,.
\end{align*}
Therefore $\gamma_M$ is mono. But it is also epi because $\gamma^M$ is. 
Therefore $\gamma_M$ is iso, and so is $\gamma^M$.
\end{proof}

\subsection{The action picture}

For a right bialgebroid $H$ over $R$ and an $H$-module algebra $M$ there are 
canonical maps
\begin{align}
\Gamma^M&:M\oR H\to\End(M_N)\,\qquad m\oR h\mapsto\{m'\mapsto m(m'\ract h)\}\\
\Gamma_M&:H\oR M\to\End(\,_NM)\,\qquad h\oR m\mapsto\{m'\mapsto (m'\ract h)m\}
\end{align}
being algebra maps from the smash products $M\mash H^\op$ and $H\mash M$, 
respectively, where in the latter case $\End(\,_NM)$ is considered with 
multiplication that arises from its natural right action on $M$. 

Note that if 
$H$ is the horizontal Hopf algebroid of a DDA and the right $H$ action arises 
from a right $V$-coaction as in Lemma \ref{lem: act-coact} then $M$ being a left 
$H^\op$-module algebra is in complete agreement with the familiar Hopf algebraic 
situation since it is $H^\op$ which is the dual of $V$. 

\begin{thm&def} \label{thm: def Gal}
Let $A$ be a distributive double algebra and $M$ a right $H$-module algebra, 
equivalently a right $V$-comodule algebra, over the horizontal, resp. vertical 
Hopf algebroid of $A$. Let $N= M^H \equiv M^\coinv{V}$. Then $N\subset M$ is 
called an $A$-Galois extension if any one of the following equivalent conditions 
hold:
\begin{tabular}{llll}
(1)&$\gamma^M$ is epi.&(2)& $\gamma_M$ is epi.\\
(3)& $\gamma^M$ is iso.&(4)& $\gamma_M$ is iso.\\
(5)& $\Gamma^M$ is iso and $M_N$ is fgp.&
(6)& $\Gamma_M$ is iso and $_NM$ is fgp.
\end{tabular}
\end{thm&def}
\begin{proof}
Equivalence of the first four conditions follows from Proposition \ref{pro: epi 
implies iso} and Lemma \ref{lem: 2 gamma}.

$(3)\Rightarrow(5)$ Considering it as a right $M$-module map, $\gamma^M$ induces 
the isomorphism (of left $M$-modules)
\[
{\gamma^M}^*\colon\Hom_{-M}(M\oT V,M)\iso\Hom_{-M}(M\oN M,M)\,.
\]
If $\chi\in\Hom_{-M}(M\oT V,M)$ then $\chi(1\oT\under)\in \Hom(V_T,M_T)$ 
because
\begin{align*}
\chi(1\oT v\star t)&=\chi(1\oT v\ci\Fi_R(t))\\
&=\chi\left(j(\Fi_R(t))\zeroT\oT v\ci j(\Fi_R(t))\oneT\right)\\
&=\chi(1\oT v)j(\Fi_R(t))\,.
\end{align*}
Thus we have a well defined map (of left $M$-modules)
\begin{align}\label{map 1}
\Hom_{-M}(M\oT V,M)&\to \Hom(V_T,M_T)\\
\chi&\mapsto \chi(1\oT\under)\notag
\end{align}
We claim that this map is an isomorphism with inverse
\[
\kappa\mapsto \{m\oT v\mapsto\kappa(v\ci S^{-1}(m\oneB))m\zeroB\}
\]
This follows from the computation 
\begin{align*}
\chi(m\oT v)&=\chi(m\zeroT\oT\Fi_T(m\oneT)\star v)
=\chi(m\zeroT\oT v\ci \Fi_L\Fi_T(m\oneT))\\
&=\chi(m\zeroT\oT v\ci S^{-1}({m\oneT}\twoB)\ci {m\oneT}\oneB)\\
&=\chi({m\zeroB}\zeroT\oT v\ci S^{-1}(m\oneB)\ci{m\zeroB}\oneT)\\
&=\chi(1\oT v\ci S^{-1}(m\oneB))\,m\zeroB
\end{align*}
on the one hand and on the other hand from $\delta_M(1)=1\oB e$.
Composing the map (\ref{map 1}) with the isomorphism
\begin{align}
\Hom(V_T,M_T)&\to M\oR H\\
\kappa&\mapsto \kappa(x^j)\oR y^j\notag
\end{align}
where $x^j\oR y^j\equiv\cop_R(e)$ is the dual basis of $\Fi_R$, we obtain the 
left vertical arrow in the diagram
\begin{equation}\label{dia: g-G 1}
\begin{CD}
\Hom_{-M}(M\oT V,M)@>{\gamma^M}^*>>\Hom_{-M}(M\oN M,M)\\
@VVV @VVV\\
M\oR H@>\Gamma^M>> \End(M_N)
\end{CD}
\end{equation}
The vertical arrow on the right is the isomorphism $\sigma\mapsto 
\sigma(\under\oN 1)$ therefore the composite along the top and right is
$\chi\mapsto \chi(\under\oT e)$. The other two compose to give 
\[
\chi\mapsto \chi(1\oT x^j)\oR y^j\mapsto \chi(1\oT x^j)(\under\ract y^j) \,.
\]
In order to see commutativity of the diagram we need a calculation.
\begin{align*}
\chi(1\oT x^j)(m\ract y^j)&=
\chi(1\oT x^j)\,m\zeroT\ud\Fi_T(m\oneT\star y^j)\\
&=\chi(m\zeroT\oT x^j\ci m\oneT\ci\Fi_R\Fi_T(m\twoT\star y^j))\\
&=\chi(m\zeroT\oT x^j\ci(m\oneT\star\Fi_T(m\twoT\star y^j))\\
&=\chi(m\zeroT\oT x^j\ci (m\oneT\star y^j)\\
&=\chi(m\zeroT\oT\Fi_L\Fi_T(m\oneT))\\
&=\chi(m\zeroT\ud\Fi_T(m\oneT)\oT e)\ =\ \chi(m\oT e)\,,
\end{align*} 
where in the fifth equality we used \cite[Equation (4.16)]{Sz: DA}.
So (\ref{dia: g-G 1}) is commutative and therefore $\Gamma^M$ is an isomorphism.

The proof of $(4)\Rightarrow(6)$ goes similarly by proving commutativity of the 
diagram
\begin{equation}\label{dia: g-G 2}
\begin{CD}
\Hom_{M-}(M\oB V,M)@>{\gamma_M}^*>>\Hom_{M-}(M\oN M,M)\\
@VVV @VVV\\
H\oR M@>\Gamma_M>> \End(\,_NM)
\end{CD}
\end{equation}
with the left hand side arrow being the isomorphism $\chi\mapsto 
x^j\oR\chi(1\oB y^j)$ and the one on the right hand side being 
$\sigma\mapsto \sigma(1\oN \under)$. 

$(5)\Rightarrow(4)$ Consider the diagram
\begin{equation} \label{dia: G-g 1}
\begin{CD}
M\oN M@>\gamma_M>>M\oB V\\
@VVV @AAA\\
\Hom_{M-}(\End(M_N),M)@>{\Gamma^M}^*>>\Hom_{M-}(M\oR H,M)
\end{CD}
\end{equation}
The lower horizontal arrow is an isomorphism since $\Gamma^M$ is. The vertical 
arrow on the left, mapping $m\oN m'$ to the homomorphism $\alpha\mapsto 
\alpha(m)m'$, is an isomorphism because $M_N$ is fgp. The other vertical arrow
is the composite of two maps,
\[
\begin{CD}
\Hom_{M-}(M\oR H,M)@>>>\Hom(\,_RH,\,_RM)@>>>M\oB V
\end{CD}
\]
where the second one is the isomorphism $\kappa\mapsto \kappa(u^k)\oB v^k$ 
with $u^k\oB v^k\equiv \cop_B(i)$ denoting the dual basis of $\Fi_B$. The first 
one, $\chi\mapsto \chi(1\oR \under)$, is obviously invertible (in contrast to 
the similar map in the $(3)\Rightarrow(5)$ part) because the left $M$-module 
structure of $M\oR H$ we need here is the trivial one. It remains to show 
commutativity of (\ref{dia: G-g 1}). So we compute the action of the lower three 
arrows,
\begin{align*}
m\oN m'&\mapsto \{\alpha\mapsto \alpha(m)m'\} \mapsto\{m''\oR h\mapsto
m''(m\ract h)m'\}\\
&\mapsto \{h\mapsto(m\ract h)m'\} \mapsto (m\ract u^k)m'\oB v^k
\end{align*}
which is indeed $\gamma_M$ if we compare the right $H$-action with the right 
$V$-coaction $\delta_M$. This proves that $\gamma_M$ is invertible.

The proof of the implication $(6)\Rightarrow(3)$ can be done similarly by using 
the diagram
\begin{equation} \label{dia: G-g 2}
\begin{CD}
M\oN M@>\gamma^M>>M\oT V\\
@VVV @AAA\\
\Hom_{-M}(\End(\,_NM),M)@>{\Gamma_M}^*>>\Hom_{-M}(H\oR M,M)
\end{CD}
\end{equation}
where on the left hand side we have the map $m\oN m'\mapsto\{\alpha\mapsto 
m\alpha(m')\}$ which is an isomorphism because $_NM$ is fgp.
\end{proof}
\begin{rmk}
The terminology "right $A$-Galois extension" where $A$ is a distributive double 
algebra does not, by any means, imply that the choice of the integral $i$ in the 
vertical Hopf algebroid $V$ plays any role. This is clear from the 
coaction picture that uses $\gamma^M$ alone. Therefore we might as well call it
"right $V$-Galois extensions" which would then be in complete agreement with the 
Hopf-Galois terminology. Saying "$A$-Galois" we try to put the 
coaction and action pictures on equal footing. For example, "bottom $A$-Galois" 
and "top $A$-Galois" extensions correspond to the left and right $H$-Galois 
extensions in the Hopf-Galois language if $H$ denotes the horizontal Hopf 
algebroid of $A$.
\end{rmk}

\subsection{Weak and strong structure theorems}

For a Frobenius Hopf algebroid $A$ let $M$ be a right $H$-module algebra and $N= 
M^H $. The category  $(\M_H)_M$ of right $M$-modules in $\M_H$ (by the 
identification $\M_H=\M^V$ being the analogue of relative Hopf modules) is 
nothing but the category of right modules over the smash product,
\begin{equation}
(\M_H)_M\ =\ \M_{H\mash M}\,.
\end{equation}
Indeed, for any action $X\oR M\to X$, $x\oR m\mapsto x\cdot m$ in $\M_H$ one has
the smash product action
\[
X\o(H\mash M)\to X,\quad x\o(h\mash m)\mapsto (x\ract h)\cdot m\,.
\]
Vice versa, any $H\mash M$-module is an $H$-module and an $M$-module and the 
$M$-action is an $H$-module map. Considering $M$ as an $N$-$H\mash M$ bimodule, 
it defines an adjoint pair $F\dashv U$ of functors 
\begin{alignat*}{2}
F&:\M_N\to\M_{H\mash M}&\qquad X&\mapsto X\oN M\\
U&:\M_{H\mash M}\to \M_N&\qquad Y&\mapsto\Hom_{H\mash M}(M,Y)
\end{alignat*}
with counit and unit
\begin{align*}
\eps_Y&:\Hom_{H\mash M}(M,Y)\oN M\to M\qquad \chi\oN m\mapsto \chi(m)\\
\eta_X&:X\to\Hom_{H\mash M}(M,X\oN M)\qquad x\mapsto\{m\mapsto x\oN m\}
\end{align*}
We note that $UY$ is isomorphic to the submodule of invariants via
\begin{equation} \label{eq: Inv Y /1}
\Hom_{H\mash M}(M,Y)\iso\Hom_H(R,Y)\iso Y^H\,.
\end{equation}
\begin{lem} \label{lem: when eta eps iso}
For any $A$-extension $N\subset M$ 
\begin{enumerate}
\item if $M_{H\mash M}$ is fgp then $\eta$ is invertible,
\item if $_NM$ is fgp and $\Gamma_M$ is invertible then $\eps$ is invertible.
\end{enumerate}
\end{lem}
\begin{proof}
(1) Apply \cite[20.10]{A-F} to the last arrow in the decomposition of $\eta_X$
\[
\begin{CD} 
X@>\sim>> X\oN M@>\sim>>X\oN\Hom_{H\mash M}(M,M)@>>>
\Hom_{H\mash M}(M,X\oN M)
\end{CD}
\]
(2) Apply \cite[20.11]{A-F} to the first arrow in the decomposition of $\eps_Y$
\[
\begin{CD} 
\Hom_{H\mash M}(M,Y)\oN M@.@.\\
@VVV@.@.\\
\Hom_{H\mash M}(\Hom_{N-}(M,M),Y)@>\Hom(\Gamma_M,Y)>>
\Hom_{H\mash M}(H\mash M,Y)@>{\sim}>>Y
\end{CD}
\]
\end{proof}
\begin{thm} \label{thm: WS}
Let $A$ be a distributive double algebra.
\begin{enumerate}
\item For an $A$-extension $N\subset M$ the following conditions are equivalent:
\begin{enumerate}
\item $\eps:FU\to\M_{H\mash M}$ is an isomorphism.
\item $N\subset M$ is $A$-Galois.
\item $_NM$ is fgp and $_NM_{H\mash M}$ is faithfully balanced.
\item $M_{H\mash M}$ is a generator.
\end{enumerate}
\item  For an $A$-Galois extension $N\subset M$ the following conditions are 
equivalent:
\begin{enumerate}
\item $F\dashv U$ is an adjoint equivalence.
\item $\eta:\,_N\M\to UF$ is an isomorphism.
\item $_NM_{H\mash M}$ is a Morita equivalence bimodule.
\item $M_{H\mash M}$ is fgp.
\item $_NM$ is a generator.
\item $_NN\subset\,_NM$ is a direct summand.
\end{enumerate}
\end{enumerate}
\end{thm}
\begin{proof}
$(1a) \Leftrightarrow(1b)$: The $\Leftarrow$ follows from Lemma \ref{lem: 
when eta eps iso} (2). As for the $\Rightarrow$ direction consider $\eps_Y$ 
for $Y=M\oT V$ which is a $H\mash M$-module via
\[
(m'\oT v)\cdot (h\mash m):=m'm\zeroT\oT(v\star h)\ci m\oneT\,.
\]
This is a well-defined action due to
\[
(m\ract h)\zeroT\oT (m\ract h)\oneT =m\zeroT\oT m\oneT\star h\,.
\]
Now consider the map
\begin{equation} \label{eq: Inv Y /2}
Y^H\to M\qquad
\sum_j\ m_j\oT w_j\mapsto \sum_j\ m_j\ud\Fi_T(w_j) 
\end{equation}
which has inverse $m\mapsto m\oT e$. As a matter of fact,
\begin{align*}
m\ud\Fi_T(e)&=m\,\eta\Fi_R\Fi_T(e)=m\eta(e)=m\\
\sum_j m_j\ud\Fi_T(w_j)\oT e&=\sum_jm_j\oT\Fi_L\Fi_T(w_j)=\sum_jm_j\oT(i\ci 
w_j)\star e\\
&=\sum_jm_j\oT(i\star e\oneR)\ci(w_j\star e\twoR)\\
&=\sum_jm_j\oT e\oneR\ci(w_j\star\Fi_B\Fi_R(e\twoR))\\
&=\sum_jm_j\oT e\oneR\ci\Fi_R(e\twoR)\ci w_j=\sum_j m_j\oT w_j\,.
\end{align*}
Composing $\eps_Y$ with the inverses of (\ref{eq: Inv Y /1}) and (\ref{eq: Inv 
Y /2}) we obtain the mapping
\[
m\oN m'\mapsto (m\oT e)\oN m'\mapsto (m\oT e)\cdot m'=m{m'}\zeroT\oT{m'}\oneT
=\gamma^M(m\oN m')
\]
Therefore $\gamma^M$ is invertible.

$(1b)\Leftrightarrow(1c)$: This is Lemma \ref{lem: ext bal} (2) together with 
the Theorem \ref{thm: def Gal} (6).  

$(1c)\Leftrightarrow(1d)$: Since $N\subset M$ is an extension, $_NM$ is 
balanced by  Lemma \ref{lem: ext bal} (1). So (c) is equivalent to that $_NM$ 
is fgp and $M_{H\mash M}$ is faithful and balanced. But these are the necessary 
and sufficient conditions for (d) by \cite[Theorem 17.8]{A-F}.

$(2a)\Leftrightarrow(2b)$: This is clear from the equivalence of (1a) and 
(1b).

$(2b)\Leftrightarrow(2c)$: Consider the composite  
\[
\begin{CD}
\Hom_N(N,X)@.\Hom_N(M\amalgo{H\mash M}\Hom(\,_NM,\,_NN),X)\\
@VV{\wr}V @A{\wr}AA\\
X@>\eta>>\Hom_{H\mash M}(M,X\oN M)
\end{CD}
\]
of natural isomorphisms where the last isomorphism exists because 
$_NM$ is fgp. By the Yoneda lemma this determines an isomorphism
\[
M\amalgo{H\mash M}\Hom(\,_NM,\,_NN)\ \to\ N\qquad\in\,_N\M
\]
which is nothing but the evaluation associated to the right dual of the bimodule 
$_NM_{H\mash M}$. Postulating the usual right $N$-module structure on 
$\Hom(\,_NM,\,_NN)$ it becomes in fact an $N$-$N$-bimodule isomorphism. 
Another hom-tensor relation for fgp $_NM$ and the isomorphism $\Gamma_M$ compose 
to give
\[
\Hom(\,_NM,\,_NN)\oN M\iso\Hom(\,_NM,\,_NM)\iso H\mash M\qquad
\in\,_{H\mash M}\M_{H\mash M}\,.
\]
Thus $\Hom(\,_NM,\,_NN)$ is the inverse equivalence of $_NM_{H\mash M}$.
It follows from Morita theory that both $_NM$ and $M_{H\mash M}$ are 
progenerators which prove that $(2c)\Rightarrow(2d)$ and 
$(2c)\Rightarrow(2e)$.

$(2d)\Leftrightarrow(2b)$ follows from Lemma \ref{lem: when eta eps iso} 
(1). 

In order to show $(2e)\Rightarrow(2f)$ we use that an $N$-module $M$
is a generator iff a finite direct sum of $M$'s contains the regular 
object as a summand, i.e., there exist $N$-module maps 
$N\rarr{\iota_k}M\rarr{\pi_k}N$ such that $\sum_k\pi_k\circ\iota_k=N$. In this
case $\{m\mapsto\sum_k\pi_k(m\,\iota_k(1))\}\in\Hom(\,_NM,\,_NN)$ splits 
the inclusion $N\subset M$. The implication $(2f)\Rightarrow(2e)$ is now 
obvious.

Finally $(2e)\Rightarrow(2d)$ follows from that $_NM_{H\mash M}$ is 
faithfully balanced by Lemma \ref{lem: ext bal}. 
\end{proof}

\subsection{An intrinsic characterization of finitary Galois extensions}

\begin{thm} \label{thm: bD2F}
For an algebra extension $N\subset M$ the following conditions are equivalent.
\begin{enumerate}
\item There is a Frobenius Hopf algebroid $V$ and a coaction of $V$ on $M$ such 
that $N\subset M$ is $V$-Galois.
\item $N\subset M$ is of depth 2 and Frobenius and $M_N$ is balanced.
\end{enumerate}
\end{thm}
\begin{proof}
$(1)\Rightarrow$ $N\subset M$ \textit{is Frobenius}: Consider the composite 
\begin{equation} \label{map 3}
\begin{CD}
M\oN M@>\gamma^M>>M\oT V@>M\o S>>M\oR H@>\Gamma^M>>\End(M_N)
\end{CD}
\end{equation}
where the middle arrow is meaningful in the double algebraic picture because $V$ 
and $H$ have the same underlying $k$-module $A$ and $S(t\star a)=S(a)\star 
\Fi_B\Fi_R(t)=\Fi_R(t)\ci a$ holds for all $a\in A$, $t\in T$, see \cite[Lemma 
5.4]{Sz: DA}. Computing the value of the map (\ref{map 3}) on $m\oN m'$ we 
obtain
\begin{align*}
m{m'}\zeroT(m''\ract S({m'}\oneT))&=
m{m'}\zeroT{m''}\zeroT\ud \Fi_T({m''}\oneT\star S({m'}\oneT))\\
&=m{m'}\zeroT{m''}\zeroT\ud \Fi_T\Fi_L({m'}\oneT\ci{m''}\oneT)\\
&=m(m'm'')\zeroT\ud\Fi_T((m'm'')\oneT\star e)\\
&=m((m'm'')\ract e)
\end{align*}
Therefore (\ref{map 3}) has the familiar form $m\oN m'\mapsto m\psi m'$ in terms 
of the $N$-$N$-bimodule map $\psi=\under \ract e$ from $M$ into $N$. Since 
(\ref{map 3}) is isomorphism it follows that $\psi$ is a Frobenius homomorphism 
with dual basis obtained from $\id_M$ by applying the inverse of (\ref{map 3}) .

$(1)\Rightarrow$ $N\subset M$ \textit{is D2}: Since $_TV$ is fgp and $\gamma^M$ 
provides an $M$-$N$-bimodule isomorphism $M\oN M\iso (_MM_N)\oT V$, it 
follows that $N\subset M$ is right D2. Similarly, the existence of the 
isomorphism $\gamma_M$ and the $_BV$ being fgp imply that $N\subset M$ is left 
D2.

$(1)\Rightarrow$ $N\subset M$ \textit{is balanced}: This follows from that every 
$V$-extension is balanced, see Lemma \ref{lem: ext bal}.

$(2)\Rightarrow(1)$: The endomorphism algebra $H^\op:=\End(\,_NM_N)$ has a 
natural structure of a Frobenius Hopf algebroid, see \cite[Subsection 8.6]{Sz: 
DA} or \cite{BSz: D2F}. Moreover, the natural action of $H^\op$ on $M$ makes it 
a left $H^\op$-module algebra and the corresponding smash product $M\mash H^\op$ 
is isomorphic to $\End(M_N)$ via $\Gamma^M$ by \cite[Corollary 4.5]{K-Sz}. So 
$N\subset M$ will be $V$-Galois, for $V$ the dual of $H^\op$, provided  
$N= M^H $. But this is equivalent to $M_N$ being balanced.
\end{proof}

Note that in the presence of the Frobenius condition left D2 is equivalent to 
right D2 and in the presence of the D2 Frobenius condition $M_N$ is balanced iff 
$_NM$ is balanced.


\section{Noncommutative scalar extensions}

The Hopf algebroid $V$ making a given algebra extension $V$-Galois is highly 
nonunique. This phenomenon can be observed already for Hopf Galois extensions. 
As Greither and Pareigis have shown \cite{Greither-Pareigis} certain separable 
field extensions can be $H$-Galois for two different Hopf algebras $H$ and $H'$.
By an appropriate extension $k\subset K$ of the scalars, however, they become
isomorphic, $K\o H\cong K\o H'$, as $K$-Hopf algebras. The $k$-Hopf algebras 
$H$, $H'$ for which such a (commutative, faithfully flat) $k$-algebra $K$ exists 
are called forms of each other \cite{Pareigis: forms}.

If we admit Hopf algebroids to appear in place of Hopf algebras then an 
interesting generalization of scalar extension is provided by the 
Brzezi\'nski-Militaru theorem \cite{Brz-Mil} constructing a Hopf algebroid 
structure on the smash product $M\mash H$ if $M$ is a braided commutative 
algebra in the Yetter-Drinfeld category $_H\YD^H$ over the Hopf algebra $H$.
As we shall see the Brzezi\'nski-Militaru theorem holds also for $H$ a 
bialgebroid or Frobenius Hopf algebroid. Since the base algebra 
of $M\mash H$ is just $M$, the braided commutative algebras (BCA's) play the 
role of (noncommutative) scalars.

If $N\subset M$ is a Galois extension for some Frobenius Hopf algebroid $H$ then 
the center $C=M^N$ of the extension is a BCA over $H$ (Corollary \ref{cor: C 
BCA}) and the scalar extension $H\mash C$ is the endomorphism Hopf algebroid 
$E$ (Proposition \ref{pro: E as smash}). Therefore all Frobenius Hopf algebroids 
$H$ for which $N\subset M$ is $H$-Galois are forms of each other.

\subsection{Braided commutative algebras}

Yetter-Drinfeld modules over bialgebroids have been introduced in 
\cite{Schauenburg: ddqg}. They form a prebraided monoidal category, the 
weak center of the category of modules over the bialgebroid.
In this subsection we adapt the weak center construction to the double algebraic 
notation and describe the (braided) center $\Z(\M_H)$ as `double' 
Yetter-Drinfeld modules $^H\YD^H_H$ with two related coactions.

For a right bialgebroid $H$ over $R$ the weak 
center $\lZ(\M_H)$ is defined as follows. The objects $\bra Z,\theta\ket$ are 
$H$-modules equipped with a natural transformation $\theta_Y:Z\oR Y\to Y\oR Z$ 
satisfying
\begin{equation} \label{eq: theta}
\theta_{X\oR Y}=(X\oR\theta_Y)\circ(\theta_X\oR Y)\quad\text{and}\quad 
\theta_R=Z
\end{equation}
where the coherence isomorphisms are not written out explicitly.
An arrow $\bra Z,\theta\ket\to\bra Z',\theta'\ket$ is an $H$-module map 
$\alpha:Z\to Z'$ such that
\begin{equation} \label{eq: arrows in lZ}
(Y\oR \alpha)\circ\theta_Y=\theta'_Y\circ(\alpha\oR Y)
\end{equation}
for all objects $Y\in\M_H$. This category has a monoidal product which is 
defined for objects by
\[
\bra Z,\theta\ket\oR\bra Z',\theta'\ket=\bra Z\oR Z',(\theta_-\oR Z')\circ
(Z\oR\theta'_-\ket
\]
and for arrows by taking the ordinary tensor product in $\M_H$. The category 
$\lZ(\M_H)$ is prebraided with
\[
\beta_{\bra Z,\theta\ket,\bra Z',\theta'\ket}=\theta_{Z'}\,.
\]

Given an object $\bra Z,\theta\ket\in\lZ(\M_H)$ one can introduce
\begin{equation}
\tau:Z\to H\oR Z\,,\qquad \tau(z):=\theta_H(z\oR i)=z\coa{-1}\oR z\coa{0}
\end{equation}
which, as being the composite
\begin{equation} \label{eq: tau}
\begin{CD}
Z@>\sim>>Z\oR R@>Z\oR\Fi_B>>Z\oR H@>\theta_H>>H\oR Z\,,
\end{CD}
\end{equation}
preserves the left $R$-module structures inherited from $\M_H$. By naturality of 
$\theta$, the $\tau$ determines $\theta_X$ for all $X$ by the formula
\begin{equation} \label{eq: theta from tau}
\theta_X(z\oR x)= x\ract z\coa{-1}\oR z\coa{0}\,.
\end{equation}
Using this formula it is easy to show that (\ref{eq: theta}) implies that 
$\tau$ is coassociative and counital, thereby making $Z$ a left 
$H$-comodule. We not only have Takeuchi's centrality property 
\begin{align} \label{eq: T for tau}
\Fi_T(r)\star z\coa{-1}\oR z\coa{0}&=\Fi_T(r)\ract z\coa{-1}\oR z\coa{0}
=\theta_H(z\oR\Fi_T(r)) \notag\\
&=\theta_H((z\oR i)\ract \Fi_T(r))=\tau(z)\ract\Fi_T(r) \notag\\
&=z\coa{-1}\oR z\coa{0}\ract \Fi_T(r)\qquad r\in R,\ z\in Z
\end{align}
but also
\begin{align} 
\Fi_B(r)\star z\coa{-1}\oR z\coa{0}&=\Fi_B(r)\ract z\coa{-1}\oR z\coa{0}
=\theta_H(z\oR\Fi_B(r)) \notag\\
&=\theta_H(z\ract \Fi_T(r)\oR i)=\tau(z\ract\Fi_T(r))=\tau(z\cdot r)\,. 
\end{align}
The latter means that the right $R$-action we could construct from the left 
$R$-action - in analogy with the left action we had in Proposition \ref{pro: 
coactions are bim} for right comodules - would be the same as the original right 
$R$-module structure inherited from (\ref{eq: tau}). In other words, 
the requirement for (\ref{eq: tau}) to be an $R$-$R$-bimodule map defines a 
right $R$-action on $H\oR Z$ which is conveyed by naturality of $\theta$ and not 
by $\theta_H$ being an arrow in $_R\M_R$.

Given a left $H$-comodule $Z$ which is also a right $H$-module
(with the same underlying $R$-$R$-bimodule structure) the condition for 
(\ref{eq: theta from tau}) to determine an $H$-module map is precisely the 
Yetter-Drinfeld condition given below.

Summarizing, one has a prebraided monoidal isomorphism $\lZ(\M_H)\cong 
\,^H\YD_H$ with the following Yetter-Drinfeld category:
\begin{defi}
For a right bialgebroid $\bra H,\star,i,R,\Fi_T,\Fi_B,\cop_R,\Fi_R\ket$ the 
category $^H\YD_H$ has objects $\bra Z,\ract,\tau\ket$ where 
\begin{enumerate}
\item $\bra Z,\ract\ket$ is a right $H$-module, hence also an 
$R$-$R$-bimodule via $r\cdot z\cdot r'=z\ract(\Fi_B(r)\star\Fi_T(r'))$.
\item $\bra Z,\tau\ket$ is a left $H$-coaction, that is to say,
\begin{enumerate}
\item $\tau:Z\to H\oR Z$ is an $R$-$R$-bimodule map in the sense of
\begin{equation} \label{eq: tau bim}
(r\cdot z\cdot r')\coa{-1}\oR(r\cdot z\cdot r')\coa{0}
=\Fi_B(r')\star z\coa{-1}\star\Fi_B(r)\oR z\coa{0}\,,
\end{equation}
\item $\tau$ is coassociative and counital,
\begin{align*}
z\coa{-1}\oR {z\coa{0}}\coa{-1}\oR{z\coa{0}}\coa{0}&=
{z\coa{-1}}\oneR\oR{z\coa{-1}}\twoR\oR z\coa{0}\\
\Fi_R(z\coa{-1})\cdot z\coa{0}&=z
\end{align*}
\item $\tau$ factorizes through $H\ex{R}Z\subset H\oR Z$, i.e., 
(\ref{eq: T for tau}) holds.
\end{enumerate}
\item The action and coaction satisfy the Yetter-Drinfeld condition
\begin{equation*}
h\twoR\star(z\ract h\oneR)\coa{-1}\oR(z\ract h\oneR)\coa{0}
=z\coa{-1}\star h\oneR\oR z\coa{0}\ract h\twoR\,.
\end{equation*}
\end{enumerate}
The arrows are the $H$-module $H$-comodule maps $Z\to Z'$. The monoidal product 
of two Yetter-Drinfeld modules $Z$ and $Z'$ is $Z\oR Z'$ equipped with 
\begin{align*}
(z\oR z')\ract h&=(z\ract h\oneR)\oR(z'\ract h\twoR)\\
(z\oR z')\coa{-1}\oR(z\oR z')\coa{0}&=
{z'}\coa{-1}\star z\coa{-1}\oR(z\coa{0}\oR {z'}\coa{0})
\end{align*}
The monoidal unit is $R$ with $r\ract h=r\star h$ and $r\coa{-1}\oR 
r\coa{0}=\Fi_B(r)\oR e$. The prebraiding is defined by
\[
\beta_{Z,Z'}:Z\oR Z'\to Z'\oR Z,\qquad z\oR z'\mapsto z'\ract z\coa{-1}\oR 
z\coa{0}\,.
\]
\end{defi}

There is a coopposite version $\rZ(\M_H)= \lZ(\M_H^\coop)=\lZ(\M_{H^\coop})$ of 
the left weak center, called the right weak center, in which an object $\bra 
Z,\bar\theta\ket$ has natural transformation $\bar\theta_Y:Y\oR Z\to Z\oR Y$ 
satisfying $\bar\theta_{X\oR Y}=(\bar\theta_X\oR Y)\circ(X\oR\bar\theta_Y)$.
This determines a right coaction
\[
\bar\tau:Z\to Z\oR H\,,\qquad z\mapsto z\coa{0}\oR z\coa{1}=\bar\theta_H(i\oR z)
\]
and is determined by this coaction,
\begin{equation} \label{eq: bartheta from bartau}
\bar\theta_Y(y\oR z)=z\coa{0}\oR y\ract z\coa{1}\,.
\end{equation}
The center $\Z(\M_H)$ is the full subcategory of $\lZ(\M_H)$ in which the 
objects $\bra Z,\theta\ket$ have invertible $\theta$. For such objects $\bra 
Z,\theta^{-1}\ket$ is an object in $\rZ(\M_H)$ in which $\bar\theta$ is 
invertible. The center is braided monoidal. In the language of 
Yetter-Drinfeld modules the objects of the center are two-sided Yetter-Drinfeld 
modules $\bra Z,\ract,\tau,\bar\tau\ket\in\,^H\YD^H_H$ in which the two 
coactions are inverse to each other, i.e.,
\begin{align}
{z\coa{0}}\coa{0}\oR z\coa{-1}\star{z\coa{0}}\coa{1}&=z\oR i
\label{eq: tau-bartau 1}\\
z\coa{1}\star{z\coa{0}}\coa{-1}\oR{z\coa{0}}\coa{0}&=i\oR z\,.
\label{eq: tau-bartau 2}
\end{align}

\begin{defi}
For a right bialgebroid $H$ the commutative monoids in $\Z(\M_H)$ are called 
BCA's (braided commutative algebras) over $H$. The commutative monoids in 
$\lZ(\M_H)$ and $\rZ(\M_H)$ are called left and right pre-BCA's over $H$, 
respectively. 
\end{defi}
Therefore a left pre-BCA consists of
an algebra $Q$ with an algebra map $\eta:R\to Q$ and
a Yetter-Drinfeld module structure $\bra Q,\ract,\tau\ket\in\,^H\YD_H$ 
such that
\begin{align}
\eta(r)\,q\,\eta(r')&=r\cdot q\cdot r'\\
(qq')\ract h&=(q\ract h\oneR)(q'\ract h\twoR)\\
1\ract h&=\eta\,\Fi_R(h)\\
(qq')\coa{-1}\oR(qq')\coa{0}&={q'}\coa{-1}\star q\coa{-1}\oR 
q\coa{0}{q'}\coa{0}\\
\eta(r)\coa{-1}\oR \eta(r)\coa{0}&=\Fi_B(r)\oR 1
\end{align}
and the prebraided commutativity
\begin{equation}
(q'\ract q\coa{-1})q\coa{0}=qq'
\end{equation}
holds. If $Q$ is a BCA then there exists also a right coaction $\bar\tau$ with 
which $\bra Q,\ract,\bar\tau\ket\in\YD^H_H$ and which is inverse to $\tau$ in 
the sense of equations (\ref{eq: tau-bartau 1}), (\ref{eq: tau-bartau 2}).

We note that the ground ring $R$ of the bialgebroid is always a BCA with the 
structure $\bra R,\mu_R,R\ket$ that comes from $R$ being the monoidal unit of 
$\lZ(\M_H)$.

\subsection{The centralizer of a Galois extension}

Interesting examples for BCA's are obtained from considering 
centralizers $M^N$ of Galois extensions. 
\begin{pro}
Let $M$ be a monoid in $\M_H$ over the right bialgebroid $H$ and let $N=M^H$.
Assume that $H_R$ is fgp and that the canonical map
$\Gamma_M:H\mash M\to\End(\,_NM)$ is an isomorphism.
Then the centralizer $M^N=\{c\in M\,|\,nc=cn,\ n\in N\}$ of the extension 
$N\subset M$ is a left pre-BCA over $H$ with $H$-module algebra structure 
inherited from $M^N\subset M$ (the Miyashita-Ulbrich action) and with left 
coaction $\tau(c):=\Gamma_M^{-1}(\lambda_M(c))$ where $\lambda_M(c)=\{m\mapsto 
cm\}$.
\end{pro}
\begin{proof}
For each $h\in H$ the action $\under\ract h$ is an $N$-$N$-bimodule map.
Therefore $M^H\subset M$ is a sub-$H$-module algebra. As such the unit 
$\eta:R\to M$ has image in $M^H$. 
Since $\Gamma_M$ is an $N$-$N$-bimodule map, it restricts to an isomorphism 
$(H\oR M)^N\iso\End(\,_NM_N)$ between the centralizers. The $H_R$ being fgp we 
have $(H\oR M)^N=H\oR M^N$. Since $\lambda(c)$ for $c\in M^N$ belongs to 
$\End(\,_NM_N)$, the $\tau$ is a map $M^N\to H\oR M^N$. The $\tau$ is uniquely 
determined by the equation
\begin{equation}  \label{eq: tau in action}
(m\ract c\coa{-1})c\coa{0}=cm\,,\qquad m\in M
\end{equation}
from which the bimodule property (\ref{eq: tau bim}) and the centrality 
(\ref{eq: T for tau}) easily follow. The calculation 
\begin{align*}
c(mm')&=((mm')\ract c\coa{-1})c\coa{0}=(m\ract {c\coa{-1}}\oneT)
(m'\ract{c\coa{-1}}\twoT)c\coa{0}\\
(cm)m'&=(m\ract c\coa{-1})c\coa{0}m'=
(m\ract c\coa{-1})(m'\ract{c\coa{0}}\coa{-1}){c\coa{0}}\coa{0}
\end{align*}
will imply coassociativity after verifying the next
\begin{lem} \label{lem: isomorphisms}
Under the assumptions of the Proposition and with the 
notations $E:=\End(\,_NM_N)$, $C:=M^N$ the maps
\begin{align}
E\oC E&\to \Hom_{N\text{-}N}(M\oN M,M) \label{eq: fork rule}\\
\alpha\oC\alpha'&\mapsto\{m\oN m'\mapsto\alpha(m)\alpha'(m')\} \notag\\
H\oR H\oR C&\to \Hom_{N\text{-}N}(M\oN M,M) \label{eq: fine fork}\\
h\oR h'\oR c&\mapsto\{m\oN m'\mapsto(m\ract h)(m'\ract h')c\}\notag
\end{align}
are isomorphisms.
\end{lem}
\begin{proof}
Using both the isomorphism $\Gamma_M$ and its restriction $H\mash C\iso E$  
we have a sequence of isomorphisms
\begin{align*}
E\oC E&\iso(H\oR C)\oC E\iso H\oR E= H\oR\Hom_{N\text{-}N}(M,M)\\
&\iso\Hom_{N\text{-}N}(M,H\oR M)\iso\Hom_{N\text{-}N}(M,\Hom_{N\text{-}}(M,M))\\
&\iso\Hom_{N\text{-}N}(M\oN M,M)
\end{align*}
The action of these isomorphisms can be computed by inserting 
$\alpha=(\under\ract h)c$ and $\alpha'=(\under\ract h')c'$:
\begin{align*}
\alpha\oC\alpha'&\mapsto(h\oR c)\oC\alpha'\mapsto h\oR c\,\alpha'(\under)\mapsto
\{m\mapsto h\oR c\,\alpha'(m)\}\\
&\mapsto\{m\mapsto\{m'\mapsto\alpha(m')\alpha'(m)\}\}\mapsto
\{m'\oN m\mapsto\alpha(m')\alpha'(m)\}
\end{align*}
This proves that (\ref{eq: fork rule}) is an isomorphism. The map in (\ref{eq: 
fine fork}) is the composite
\begin{equation*}
\begin{CD}
H\oR H\oR C@.@.\Hom_{N\text{-}N}(M\oN M,M)\\
@V{H\oR\Gamma}VV @.  @AA{\cong}A\\
H\oR E@>\cong>>H\oR C\oC E@>\Gamma\oC E>>E\oC E
\end{CD}
\end{equation*}
of isomorphisms. 
\end{proof}
Returning to the proof of the Proposition counitality of $\tau$
can be seen as
\[
\Fi_R(c\coa{-1})\cdot c\coa{0}=(1\ract c\coa{-1})c\coa{0}=c1=c\,.
\]
As for the Yetter-Drinfeld compatibility condition 
it suffices to verify the equality
\begin{align*}
(m\ract c\coa{-1}\star h\oneT)(c\coa{0}\ract h\twoT)&=\left((m\ract 
c\coa{-1})c\coa{0}\right)\ract h=(cm)\ract h\\
=(c\ract h\oneT)(m\ract h\twoT)&=\left(m\ract h\twoT\star(c\ract 
h\oneT)\coa{-1}\right)(c\ract h\oneT)\coa{0}
\end{align*}
In order to see compatibility of $\tau$ with multiplication and unit in $C$ it 
suffices to check
\[
c'cm=c'(m\ract c\coa{-1})c\coa{0}=
(m\ract c\coa{-1}\star {c'}\coa{-1}){c'}\coa{0}c\coa{0}\,.
\]
Finally, braided commutativity $(c'\ract c\coa{-1})c\coa{0}=cc'$
follows from the more general relation (\ref{eq: tau in action}).
\end{proof}

\begin{cor} \label{cor: C BCA}
If $N\subset M$ is a right $A$-Galois extension for a distributive double 
algebra $A$ then $M^N$ is a BCA over the horizontal Hopf algebroid $H$.
\end{cor}
\begin{proof}
It suffices to prove that the prebraiding is invertible. Define the right 
coaction $\bar\tau(c):=(\Gamma^M)^{-1}(\rho_M(c))$ where $\rho_M$ is 
right multiplication on $M$. This is equivalent to $\bar\tau(c)=c\coa{0}\oR 
c\coa{1}$ satisfying
\begin{equation} \label{eq: bartau in action}
c\coa{0}(m\ract c\coa{1})=mc\,,\qquad m\in M\,.
\end{equation}
Applying (\ref{eq: tau in action}) to (\ref{eq: bartau in action}) we obtain
\[
(m\ract i)c=mc=(m\ract c\coa{-1}\star{c\coa{0}}\coa{-1})\,{c\coa{0}}\coa{0}
\]
from which equation (\ref{eq: tau-bartau 2}) follows. Equation (\ref{eq: 
tau-bartau 1}) can be seen similarly.
\end{proof}
Notice that this proof does not use very much from the Hopf algebroid structure.
Therefore the Corollary holds true for any right bialgebroid for which both 
$H_R$ and $_RH$ are fgp and for all extensions for which both $\Gamma^M$ and 
$\Gamma_M$ are invertible.

\subsection{Extensions by BCA's}

For any $H$-module algebra $Q$ over the right bialgebroid $H$ the category
$\M_{H\mash Q}$ of modules over the smash product can be identified with the
category of (internal) $Q$-modules $(\M_H)_Q$ in $\M_H$.

If $Q$ is also a pre-BCA then every right $Q$-module in $\M_H$ is also a left 
$Q$-module by pre-braided commutativity. This defines an embedding of categories
\begin{equation} \label{eq: BM 1}
\M_{H\mash Q}=(\M_H)_Q\hookrightarrow \,_Q(\M_H)_Q
\end{equation}
into the monoidal category of internal $Q$-$Q$-bimodules. Since the 
$Q$-$Q$-bimodule tensor product of diagonal bimodules $X,Y\in(\M_H)_Q$ is again 
diagonal due to one of the hexagons, this embedding is actually strong monoidal.
Composing (\ref{eq: BM 1}) with the strong monoidal forgetful functor 
$_Q(\M_H)_Q\to \,_Q\M_Q$ we obtain a strong monoidal functor
\begin{equation} \label{eq: BM 2}
\M_{H\mash Q}=(\M_H)_Q\to \,_Q\M_Q\,.
\end{equation}
This functor is precisely the forgetful functor associated to the algebra map
\begin{equation} \label{eq: forget along}
Q^\op\o Q\to H\mash Q\,, \qquad q\o q'\mapsto q\coa{-1}\mash q\coa{0}q'
\end{equation}
therefore, by a theorem of Schauenburg \cite{Schauenburg: bnrsthb}, there is
a unique bialgebroid structure on $H\mash Q$ such that the given monoidal 
structure of $\M_{H\mash Q}$ is that of the module category of a bialgebroid. 
This is the Brzezi\'nski-Militaru Theorem in disguise. More precisely this is 
the "only if" part of \cite[Theorem 4.1]{Brz-Mil} generalized to bialgebroids 
$H$.\begin{thm}
Let $H$ be a right bialgebroid over $R$ and let $Q$ be a left pre-BCA over $H$. 
Then the smash product $G:=H\mash Q$ is a right bialgebroid over $Q$ with 
structure maps
\begin{align}
s_G(q)&=i\mash q \label{s_G}\\ 
t_G(q)&=q\coa{-1}\mash q\coa{0}\label{t_G}\\
\cop_G(h\mash q)&=(h\oneR\mash 1)\oQ(h\twoR\mash q)\label{cop_G}\\
\eps_G(h\mash q)&=\eta(\eps_H(h))q\label{eps_G}
\end{align}
where $\eta:R\to Q$ is the unit of $Q$. Moreover, $h\mapsto h\mash 1$ is a
bialgebroid map $\iota:H\to G$.

If $H$ is a Frobenius Hopf algebroid with Frobenius integral $e$ then $G$ is 
also a Frobenius Hopf algebroid with $e_G=\iota(e)$ a Frobenius integral.
\end{thm}
\begin{proof}
The observation made before the 
formulation of the Theorem, in particular equation (\ref{eq: forget along}) 
implies the formulae for $s_G$ and $t_G$. 
In order to obtain the expressions for $\cop_G$ and $\eps_G$ at once,
and also to prove the Frobenius Hopf algebroid case, the next Proposition, 
however simple, is very useful.
\begin{pro} \label{pro: tensor Q}
If $H$ is a right bialgebroid over $R$ and $Q$ is a left pre-BCA over $H$ then
the functor $\under\oR Q:\M_H\to(\M_H)_Q$ is strong monoidal.
\end{pro}
\begin{proof}
The natural transformation
\begin{align*}
(Y\oR Q)\oQ(Y'\oR Q)&\to (Y\oR Y')\oR Q\\
(y\oR q)\oQ(y'\oR q')&\mapsto(y\oR y'\ract q\coa{-1})\oR q\coa{0}q'
\end{align*}
has inverse $(y\oR y')\oR q\mapsto (y\oR 1)\oQ(y'\oR q)$. The $H\mash Q$-module 
map
\[
Q\to R\oR Q\,,\qquad q\mapsto e\oR q
\]
is the unit part of the monoidal structure and is obviously invertible.
\end{proof}
Continuing the proof of the Theorem we take the comonoid $\bra 
H,\cop_H,\eps_H\ket$ in $\M_H$ and apply the strong monoidal functor $\under\oR 
Q$. It is easy to check that the result is precisely $\bra G,\cop_G,\eps_G\ket$ 
which is then necessarily a comonoid in $\M_G$. This comonoid is obviously 
strong \cite{Sz: MonMor} proving that $\bra G,Q,s_G,t_G,\cop_G,\eps_G\ket $ is a 
bialgebroid. It is straightforward to verify that the pair $\bra \iota,\eta\ket$ 
satisfies the four axioms \cite{Sz: Strasbourg,Sz: MonMor} for a bialgebroid map 
$H\to G$.


If $H$ is a Frobenius Hopf algebroid then it has a distributive double algebra 
structure \cite{Sz: DA}. Therefore we may assume that $H$ is the horizontal Hopf 
algebroid of $\bra A,\ci,e,\star,i\ket$. Then $\bra 
H,\cop_R,\Fi_R,\ci,R\hookrightarrow H\ket$ is a Frobenius algebra in $\M_H$, so 
it is mapped by the strong monoidal functor of Proposition \ref{pro: tensor Q} 
to a Frobenius algebra in $\M_G$. The comonoid part of this Frobenius algebra 
has already been determined to be $\bra G,\cop_G,\eps_G\ket$. The monoid part 
will provide a convolution product with unit on $G$ which, together with the 
smash product algebra structure, will make $G$ a distributive double algebra. 
This convolution product (vertical multiplication) is obtained as the composite
\[
(h\mash q)\oQ(h'\mash q')\mapsto (h\oR h'\star q\coa{-1})\oR q\coa{0}q'
\mapsto h\ci(h'\star q\coa{-1})\mash q\coa{0} q'
\]
and its unit element $e_G$ is the image of $1\in Q$ under the map
\[
Q\iso R\oR Q\to H\mash Q \,.
\]
So $e_G=e\mash 1$ is a two-sided Frobenius integral in $G$.
\end{proof}

\begin{rmk}
The construction of a vertical multiplication on $H\mash Q$ suggests the new 
interpretation of the smash product as a double algebraic one. If $\bra 
A,\ci,e,\star,i\ket$ is a DDA and $Q$ is a BCA over the bialgebroid 
$H$ over $R$ then there is a smash product double algebra $A\mash Q$ with
\begin{itemize}
\item underlying $k$-module $A\oR Q$,
\item horizontal multiplication
$
(a\mash q)\star(a'\mash q')=a\star {a'}\oneR\mash(q\ract {a'}\twoR)q'\,,
$
\item horizontal unit $i\mash 1$,
\item vertical multiplication
$
(a\mash q)\ci(a'\mash q')=a\ci(a'\star q\coa{-1})\mash q\coa{0}q'\,,
$
\item and vertical unit $e\mash 1$.
\end{itemize}
\end{rmk}

As a biproduct of the double algebraic picture we obtain the following result.
\begin{pro} \label{pro: YD braided}
For Frobenius Hopf algebroids $H$ the prebraiding of the left weak center 
$\lZ(\M_H)$ is a braiding. Therefore $\lZ(\M_H)=\Z(\M_H)=\rZ(\M_H)$ and every 
pre-BCA is a BCA over $H$.
\end{pro}
\begin{proof}
We claim that the inverse braiding encoded in the right coaction $\bar\tau$ by 
(\ref{eq: bartheta from bartau}) is given by
\begin{equation} \label{eq: inverse tau}
q\coa{0}\oR q\coa{1}=\eta_Q\Fi_R\Fi_T(x^j\star q\coa{-1})q\coa{0}\oR y^j\,.
\end{equation}
The proof is motivated by the double algebraic structure on $H\mash Q$
given in the above Remark but we do not use that the given structure maps 
satisfy the axioms of a DDA. Let us compute the would-be $\Fi_R$ of $H\mash Q$.
It is
\[
\Phi_R(h\mash q):=(e\mash 1)\star(h\mash q)=e\mash \eta_Q\Fi_R(h)q\,.
\]
One conjectures $(x^j\mash 1)\oQ (y^j\mash 1)$ to be its dual basis.
Instead of proving that we prove its special case
\begin{align*}
\Phi_R((i\mash q)\ci(x^j\mash 1))\ci (y^j\mash 1)
&=(e\mash\eta_Q\Fi_R\Fi_T(x^j\star q\coa{-1})q\coa{0})\ci(y^j\mash 1)\\
&=y^j\star(\Fi_R\Fi_T(x^j\star q\coa{-2})\ci q\coa{-1})\mash q\coa{0}\\
&=y^j\star q\coa{-1}\star \Fi_B\Fi_R\Fi_T(x^j\star q\coa{-2})\mash q\coa{0}\\
&=\Fi_R(i\ci\Fi_R(q\coa{-1})\oneR)\ci\Fi_R(q\coa{-1})\twoR\mash q\coa{0}\\
&=i\ci \Fi_R(q\coa{-1})\mash q\coa{0}\\
&=i\mash q\,.
\end{align*}
Comparing the first row with the Ansatz (\ref{eq: inverse tau}) and then using 
the vertical multiplication of $H\mash Q$ we arrive at
\begin{align*}
i\mash q&=(e\mash q\coa{0})\ci(q\coa{1}\mash 1)\\
&=q\coa{1}\star{q\coa{0}}\coa{-1}\mash{q\coa{0}}\coa{0}
\end{align*}
which is equation (\ref{eq: tau-bartau 2}). The verification of  
(\ref{eq: tau-bartau 1}) is a bit longer,
\begin{align*}
{q\coa{0}}\coa{0}\oR q\coa{-1}\star{q\coa{0}}\coa{1}
&=\eta_Q\Fi_R\Fi_T(x^j\star q\coa{-1})q\coa{0}\oR q\coa{-2}\star y^j\\
&=\eta_Q\Fi_R\Fi_T(S^{-1}(q\coa{-2})\star x^j\star q\coa{-1})q\coa{0}\oR y^j\\
&=\eta_Q\Fi_R\Fi_T(S^{-1}(x^k)\star x^j\star(y^k\ci q\coa{-1}))q\coa{0}\oR y^j\\
&=\eta_Q\Fi_R\Fi_T((S^{-1}(y^k)\ci(S^{-1}(x^k)\star x^j))\star q\coa{-1})
q\coa{0}\oR y^j\\
&=\eta_Q\Fi_R\Fi_T((x_k\ci(y_k\star x^j))\star q\coa{-1})q\coa{0}\oR y^j\\
&=\eta_Q\Fi_R\Fi_T(\Fi_R\Fi_T(x^j)\star q\coa{-1})q\coa{0}\oR y^j\\
&=\eta_Q\Fi_R(\Fi_T(x^j)\star q\coa{-1})q\coa{0}\oR y^j\\
&=\eta_Q\Fi_R(q\coa{-1})q\coa{0}\ract \Fi_T(x^j)\oR y^j\\
&=q\oR\Fi_R\Fi_T(x^j)\ci y^j\\
&=q\oR i\,.
\end{align*}
\end{proof}

Extensions of quantum groupoids by BCA's are transitive in the following 
sense.
\begin{pro}
If $Q$ is a BCA over $H$ and $P$ is a BCA over $H\mash
Q$ then $P$ is a BCA over $H$, too. Furthermore, 
$ (H\mashed{R} Q)\mashed{Q} P \cong H\mashed{R} P$.
\end{pro}
\begin{proof} (Sketch)
The notation $P$ comprises the data $\bra\bra P,\theta^P\ket,\mu_P,\eta_P\ket$ 
of a commutative monoid in $\lZ(\M_{H\mash Q})$. We use analogous notations for 
$Q$. In order to obtain a monoid $\bra\bra P,\theta\ket,\mu,\eta\ket$ in 
$\lZ(\M_H)$ we define
\begin{align}
\theta_X&:
\begin{CD}
P\oR X\iso P\oQ Q\oR X@>P\o\theta^Q_X>>P\oQ(X\oR Q)@>\theta^P_{X\o Q}>>
X\oR Q\oQ P\iso X\oR P
\end{CD}\\
\mu&:P\oR P\to P\oQ P\rarr{\mu_P}P\\
\eta&:R\rarr{\eta_Q}Q\rarr{\eta_P}P
\end{align}
For the $\theta_X$ to be a well-defined arrow in $\M_H$ notice that $\theta^Q_X$ 
belongs to $_Q\M_H$ due to that $\mu_Q$ satisfies (\ref{eq: arrows in lZ}) and 
is braided commutative.The pair $\bra P,\theta\ket$ will be established as an 
object in the left weak center if we can show that (\ref{eq: theta}) is 
satisfied. This is rather weary but straightforward utilizing that (\ref{eq: 
theta}) is satisfied by both $\theta^Q$ and $\theta^P$. It is clear that the 
triple $\bra P,\mu,\eta\ket$ is a monoid in $\M_H$, but we need to show that 
$\mu$ and $\eta$ belong to $\lZ(\M_H)$. This is equivalent to checking  
(\ref{eq: arrows in lZ}), i.e., the equations
\begin{align}
(X\oR\mu)\circ(\theta_X\oR P)\circ(P\oR \theta_X)&=\theta_X\circ(\mu\oR X)\\
(X\oR \eta)\circ\theta^R_X&=\theta_X\circ(\eta\oR X)\,.
\end{align}

\end{proof}

A sort of converse to the previous proposition is the next proposition
which we state without proof.
\begin{pro}
If $\eta:Q\to P$ is a monoid morphism in $\Z(\M_H)$ between commutative 
monoids (i.e., BCA's over $H$) then there is a unique BCA structure on 
$P$ over $H \mash Q$ the unit of which is $\eta$. 
\end{pro}


\begin{pro} \label{pro: E as smash}
Let $N\subset M$ be a Galois extension over the Frobenius Hopf algebroid $H$. 
Then the restriction of the Galois map $\Gamma_M$ provides an isomorphism of 
Hopf algebroids $H\mash C\cong E$ where $E$ is the endomorphism Hopf algebroid 
of the extension.
\end{pro}
\begin{proof}
The structure maps (\ref{s_G}), (\ref{t_G}), (\ref{cop_G}) and (\ref{eps_G}) of 
the smash product are mapped by $\Gamma_M$ to 
\begin{alignat}{2}
s_E&:C\to E&\qquad c&\mapsto \{m\mapsto mc\}\\
t_E&:C^\op\to E&\qquad c&\mapsto \{m\mapsto cm\}\\
\cop_E&:E\to E\oC E\quad &\text{such that } 
&\alpha\oneR(m)\alpha\twoR(m')=\alpha(mm')\\
\eps_E&:E\to C&\qquad \alpha&\mapsto \alpha(1)
\end{alignat}
respectively, where note that multiplicativity of $\cop_E$ uniquely fixes it by 
Lemma \ref{lem: isomorphisms}, (\ref{eq: fork rule}).
Now it is easy to check that $\Gamma_M:H\mash C\to E$ satisfies the axioms 
of bialgebroid maps.
\end{proof}

\section{Contravariant fiber functors} \label{s: fiber}

In this section we study functors from the module category of a Hopf algebroid 
$A$ that correspond to $A$-Galois extensions of a given algebra $N$.
In this sense we study generalizations of Ulbrich's Theorem \cite{Ulbrich} 
relating Hopf-Galois extensions to fiber functors.
Technically speaking, however, the functors we study here are very different 
from the usual fiber functors. They are contravariant hom-functors 
$\Hom_H(\under,M)$ from $\M_H$ to $_N\M_N$. So they are colimit preserving but 
rarely faithful and exact. Still they have some properties that are worthy of
discussion. As a preparation we prove 
\begin{lem} \label{lem: fgp mon}
For $H$ a Frobenius Hopf algebroid the full subcategory $\M_H^\fgp$ of $\M_H$ 
the objects of which are the finitely generated projective $H$-modules is a 
monoidal subcategory.
\end{lem}
\begin{proof}
It suffices to show that $H\oR H$, the tensor square of the regular object in 
$\M_H$ is a fgp module. This in turn follows from the existence of the 
isomorphism \cite[(4.1)]{Sz: DA}
\[
\Gamma_{RB}:H\oR H\iso A\oB H\,,\qquad a\oR a'\mapsto a\oneB\oB a\twoB\ci a'\,.
\]
which happens to be an $H$-module map,
\begin{align*}
\Gamma_{RB}(a\ract h\oneR\oR a'\ract {h'}\twoR)
&=a\star h\oneR\star u_k\oB v_k\ci(a'\star h\twoR)\\
&=a\star u_k\oB (v_k\ci a')\star h
\end{align*}
thanks to right distributivity in $A$. Since $A_B$ is fgp, the statement is 
proven.
\end{proof}

\begin{thm} \label{thm: U}
Let $N$ be an algebra and $A$ a distributive double algebra. As usual, $H$ 
denotes the horizontal Hopf algebroid of $A$.
\begin{enumerate}
\item The mappings
\[
M\mapsto F= \Hom_H(\under,M)\quad\text{respectively}\quad F\mapsto M=F(H_H)
\]
provide mutually inverse category equivalences between the following two 
categories.
\begin{itemize}
\item The category of $H$-module algebras $M$ equipped with an algebra map 
$\hat\eta:N\to  M^H \subset M$. The arrows 
from $\bra M,\hat\eta\ket$ to $\bra M',\hat\eta'\ket$ are the $H$-module algebra 
maps $\alpha:M\to M'$ for which $\alpha\circ\hat\eta=\hat\eta'$.
\item The category of opmonoidal functors 
$F:\M_H^\fgp\to\,_N\M_N^\op$ as objects and 
monoidal natural transformations as arrows.
\end{itemize}
\item $M$ is an $A$-extension of $N$ iff $F$ is normal opmonoidal.
\item $M$ is an $A$-Galois extension of $N$ iff $F$ is strong (op)monoidal.
\item The full subcategories $A$-$\Gal(N)$ and $\F(\M_H^\fgp,\,_N\M_N^\op)$ of 
those in (1) selected by the conditions of (3), respectively, are groupoids.
\end{enumerate}
In this way (1) and (3) establish a category equivalence $A$-$\Gal(N)\sim 
\F(\M_H^\fgp,\,_N\M_N^\op)$ between $A$-Galois extensions of $N$ and strong 
monoidal  functors $\M_H^\fgp\to \,_N\M_N^\op$.
\end{thm}
\begin{proof} 
The construction of the functor $M\mapsto F$ goes as follows. Given $\bra 
M,\hat\eta\ket$ the $M$ is an $N^e$-$H$-bimodule so $\Hom_H(\under,M)$ is a 
contravariant functor from $H$-modules to $_N\M_N$.The monoid structure $\bra 
M,\mu,\eta\ket$ defines an opmonoidal structure on this functor\footnote{The 
arrows between $N$-$N$-bimodules are always considered in $_N\M_N$ and never in 
$_N\M_N^\op$.} 
\begin{align*}
\Hom_H(Y,M)\oN\Hom_H(Y',M)&\to\Hom_H(Y\oR Y',M)\quad 
\xi\oN\xi'\mapsto\mu\circ(\xi\oR\xi')\\
N&\to\Hom_H(R,M)\qquad n\mapsto\{r\mapsto n\eta(r)\}
\end{align*}
An arrow $\alpha$ is mapped to the monoidal natural transformation 
$\Hom_H(\under,\alpha)$. 

Now we construct the functor $F\mapsto M$.
Any opmonoidal functor $\bra F,F^2,F^0\ket:\M_H^\fgp\to\,_N\M_N^\op$
maps comonoids to comonoids. Therefore it maps $\bra H_H,\cop_R,\Fi_R\ket$ to a 
monoid $M:=F(H_H)$ in $_N\M_N$. The unit of this monoid is the 
composite
\[
\hat\eta:=\quad 
\begin{CD}
N@>F^0>> FR@>F(\Fi_R)>> M
\end{CD}
\]
which becomes a $k$-algebra map by prolongation of the multiplication 
\[
\hat\mu:=\quad
\begin{CD}
M\oN M@>F^{H,H}>>F(H\oR H)@>F(\cop_R)>> M
\end{CD}
\]
to a $k$-algebra multiplication. The $M$ also inherits a right $H$-module 
structure from left multiplication $\lambda_h=h\star\under$ via 
$H\rarr{\lambda}\End(H_H)\rarr{F}\End(\,_NM_N)$. We have
\begin{align*}
F\lambda_h\circ\hat\eta&=F(\Fi_R\circ\lambda_h)\circ 
F^0=F(\Fi_R\circ\lambda_{\Fi_T\Fi_R(h)})\circ F^0\\
&=F\lambda_{\Fi_T\Fi_R(h)}\circ\hat\eta
\end{align*}
implying that $\hat\eta$ factors uniquely through the inclusion 
$ M^H \subset M$. 
For each $h$ the action $F\lambda_h$ is an $N$-$N$ bimodule 
map which makes $M$ an $N^e$-$H$-bimodule. By means of the isomorphism 
$\iota:m\mapsto \{h\mapsto m\ract h\}$ the monoid 
$\bra\,_NM_N,\hat\mu,\hat\eta\ket$ becomes isomorphic to the convolution monoid 
$\Hom_H(H,M)$ associated to an $H$-module algebra $\bra M,\mu,\eta\ket$ 
structure on $M$. Of course, the monoid $\bra M,\mu,\eta\ket$ arises from the
$k$-algebra structure of $M$ just as the monoid $\bra M,\hat\mu,\hat\eta\ket$ 
does.
\begin{equation} \label{dia: mul-mul}
\begin{CD}
FH\o FH@>\iota\o\iota>>\Hom_H(H,M)\o\Hom_H(H,M)\\
@VVV @VVV\\
FH\oN FH@>\iota\oN\iota>>\Hom_H(H,M)\oN\Hom_H(H,M)\\
@V{F^{H,H}}VV @VV{\mu\circ(\under\oR\under)}V\\
F(H\oR H)@.\Hom_H(H\oR H,M)\\
@V{F\cop_R}VV @VV{\under\circ\cop_R}V\\
FH@>\iota>>\Hom_H(H,M)
\end{CD}
\end{equation}
This yields the object map of the functor $F\mapsto M$. As 
for the arrow map take any monoidal natural transformation $\nu:F\to F'$ and 
define $\alpha:=\nu_H:M\to M'$. Then by the multiplicativity constraint 
for $F$ the $\alpha$ is an $H$-module algebra morphism and the unit constraint 
implies that
\[
\begin{CD}
N@>F^0>>FR@>F\Fi_R>>M\\
@| @VV{\nu_R}V @VV{\alpha}V\\
N@>{F'}^0>>F'R@>F'\Fi_R>>M'
\end{CD}
\]
is commutative, i.e., $\alpha\circ\hat\eta=\hat\eta'$. 

Now we construct a natural isomorphism $\nu$ from the identity functor $F\mapsto 
F$ to the composite $F\mapsto M\mapsto F$. Choosing a direct 
summand diagram $Y\rarr{\pi_k}H\rarr{\sigma_k}Y$ for each fgp $H$-module $Y$
the isomorphism $\iota:M\to \Hom_H(H,M)$ for $M=FH$ extends to a natural 
isomorphism $\nu:F\to\Hom_H(\under,M)$ by
\[
\begin{CD}
FY@>F\sigma_k>>FH@>F\pi_k>>FY\\
@V{\nu_Y}VV @VV{\iota}V @VV{\nu_Y}V\\
\Hom_H(Y,M)@>\under\circ\sigma_k>>\Hom_H(H,M)@>\under\circ\pi_k>>
\Hom_H(Y,M)
\end{CD}
\]
This natural isomorphism will then be 
automatically monoidal due to the interplay between the multiplications 
$\hat\mu$ and $\mu$ seen on the diagram (\ref{dia: mul-mul}) .

The natural isomorphism from the identity functor $M\mapsto M$ to the composite 
$M\mapsto F\mapsto M$ is just $\iota:M\to\Hom_H(H,M)$ viewed as a map of monoids 
in $_N\M_N$. So in particular $\iota\circ\hat\eta$ is equal to $F\Fi_R\circ F^0$ 
for the opmonoidal functor $F=\Hom_H(\under,M)$.
This completes the proof of the equivalence in (1).

By Lemma \ref{lem: R->M} the unique arrow 
$N\to  M^H $ factorizing $\hat\eta$ is an isomorphism iff $F^0$ is an 
isomorphism, i.e., iff $F$ is normal. This proves (2).

Strong (op)monoidality of $F$ is equivalent to invertibility of $F^0$ and 
$F^{H,H}$. By the natural isomorphism $F\cong\Hom_H(H,M)$ the latter is
equivalent to invertibility of the left vertical arrow in the next diagram.
\[
\begin{CD}
\Hom_H(H,M)\oN \Hom_H(H,M)@<\sim<< M\oN M\\
@V{\mu\circ(\under\oR\under)}VV @VV{\gamma^M}V\\
\Hom_H(H\oR H,M)@>\sim>>M\oT A
\end{CD}
\]
where the lower horizontal arrow is given by a composition of isomorphisms
\[
\begin{CD}
\Hom_H(H\oR H,M)@>\sim>>\Hom_H(H,M\oR H)@>\sim>> M\oR H
@>M\o S^{-1}>> M\oT V
\end{CD}
\]
performing the mappings
\[
\chi\mapsto \chi(\under\oR x^j)\oR y^j\mapsto \chi(i\oR x^j)\oR y^j\mapsto 
\chi(i\oR u^j)\oT v^j\,.
\]
Commutativity of the diagram now follows from the simple calculation
\[
(m\ract i)(m'\ract u^k)\oT v^k=m{m'}\zeroT\oT {m'}\oneT=\gamma^M(m\oN m')
\]
Therefore $\gamma^M$ is invertible iff $F^{H,H}$ is invertible. Adding the 
condition that $N\subset M$ is an $A$-extension we obtain (3).  

Since $H_H$ is a Frobenius algebra, it is a selfdual object in $\M_H$. Therefore 
any monoidal natural transformation between strong monoidal functors from 
$\M_H^\fgp$ is invertible at $H_H$ \cite{Saavedra-Rivano} and therefore it is 
invertible everywhere. This proves (4).
\end{proof}

\begin{cor}
The map $M\mapsto \Hom_H(\under,M)$ is a category equivalence 
between the category $A$-$\Gal(N)$ of $A$-Galois extensions of $N$ and
the category $\F(\M_H,\,_N\M_N^\op)$ of colimit preserving opmonoidal functors
the restrictions of which to $\M_H^\fgp$ is strong (op)monoidal.
\end{cor}
\begin{proof}
If $F:\M_H\to \,_N\M_N^\op$ is colimit preserving then the corresponding 
$F^\op:\M_H^\op\to\,_N\M_N$ is limit preserving and $H_H$ is a cogenerator for 
$\M_H^\op$. The conditions for the special adjoint functor theorem 
\cite{MacLane} hold, so $F^\op$ has a left adjoint. It follows that $F^\op$ is a 
hom-functor, $F\cong\Hom_{\M_H^\op}(M\under)$, i.e., $F\cong\Hom_H(\under,M)$.
Now Theorem \ref{thm: U} implies that $F$ has a strong restriction to the fgp 
modules precisely when $N\subset M$ is $A$-Galois. Vice versa, every Galois 
extension $M$ gives rise to a colimit preserving opmonoidal functor 
$\Hom_H(\under, M)$ the resriction of which to $\M_H^\fgp$ is strong.
\end{proof}

As an application of the strong monoidal functor $\Hom(\under,M)$ we present 
here another characterization of Galois extensions over DDA's.
In order to understand the terminology "left distributivity"
let us look at multiplication of $M$ as a vertical one and the right $H$-action 
$\ract$ as a partially defined horizontal multiplication between $M$ and $H$.
\begin{pro}
Let $A$ be a DDA and $M$ be a right $A$-module algebra with $N=M^H$.
Then $N\subset M$ is $A$-Galois if and only if $\psi=\under\ract e:M\to N$ is a 
Frobenius homomorphism and the "left distributivity" rule
\[
m\ract (a\ci a')=(m\oneL \ract a)(m\twoL\ract a')
\]
holds for all $m\in M$ and $a,a'\in A$. Here $m\oneL\oN m\twoL$ is the coproduct
associated to the Frobenius structure on $N\subset M$ defined by $\psi$.
\end{pro} 
Note that "right distributivity" $(mm')\ract a=(m\ract a\oneR)(m'\ract a\twoR)$ 
holds for all right module algebras. Note also that $\ci$ for $H$ plays the role
of convolution product while the ordinary product is $\star$.
\begin{proof}
\textit{Necessity:} Consider the contravariant functor $\Hom_H(\under,M):\M_H\to 
\,_N\M_N$. It is strong monoidal, so maps monoids to comonoids, comonoids to 
monoids, and Frobenius algebras to Frobenius algebras. Therefore it maps 
$\bra A,\cop_R,\Fi_R,\mu_V,R\hookrightarrow A\ket$ to some
Frobenius algebra structure on $\Hom_H(A,M)\cong M\in\,_N\M_N$. Since a 
Frobenius algebra structure in $_N\M_N$ is uniquely determined by the algebra 
structure and by the Frobenius homomorphism, the counit, it is sufficient to 
check that the image of $\bra A,\cop_R,\Fi_R\ket$ is the convolution algebra 
$\Hom_H(A,M)$ and the image of the unit $R\hookrightarrow A$ is $\psi$. 
Then the coproduct must have the form
\[
\cop_M(m)\equiv m\oneL\oN m\twoL=\sum_i me_i\oN f_i
\]
where $\sum_i e_i\oN f_i$ is the dual basis of $\psi$. This means that the
composite
\begin{equation}
\begin{CD}
\Hom_H(A,M)@>\Hom(\mu_V,M)>>\Hom_H(H\oR H,M)\\
@. @V{\wr}V{[\mu\circ(\under\oR\under)]^{-1}}V\\
@.\Hom_H(A,M)\oN\Hom_H(A,M)
\end{CD}
\end{equation}
must be the map
\[
(m\ract\under)\mapsto(m\oneL\ract\under)\oN(m\twoL\ract\under)
\]
Applying $\mu\circ(\under\oR\under)$ we obtain left distributivity.

\textit{Sufficiency:} Consider the map $M\oT A\to M\oN M$ defined by $m\oT 
a\mapsto m\oneL\oN m\twoL\ract a$. Then
\begin{align*}
\gamma^M(m\oneL\oN m\twoL\ract a)&=m\oneL (m\twoL\ract a)\zeroT\oT 
(m\twoL\ract a)\oneT\\
&=(m\oneL\ract i) (m\twoL\ract u^k)\oT v^k\star a=(m\ract (i\ci u^k))\oT 
v^k\star a\\
&=m\oT \Fi_T(u^k)\star v^k\star a=m\oT a
\end{align*}
proves that $\gamma^M$ is epi.
\end{proof}

\section{A monoidal duality}

Given a right bialgebroid $H$ over $R$, an $H$-module algebra $M$ and an 
algebra map $N\to M^H $  we can 
look for a duality between - full subcategories of - $\M_H$ and $_N\M_N$ in the 
following form. The $M$ being an $N^e$-$H$-bimodule, it determines two functors
\begin{alignat}{2}
J&:\,_N\M_N^\op&\to\M_H\qquad X&\mapsto\Hom_{N^e}(X,M)\\
K&:\M_H&\to \,_N\M_N^\op\qquad Y&\mapsto\Hom_H(Y,M),
\end{alignat}
the $M$-dual functors,
that are in adjunction $K\adj J$. The counit and unit of the adjunction  
are just the natural homomorphism to the double dual,
\begin{align*}
\sigma_X&:X\to\Hom_H(\Hom_{N^e}(X,M),M)\qquad\in \,_N\M_N\\
\sigma_Y&:Y\to\Hom_{N^e}(\Hom_H(Y,M),M)\qquad\in\ \ \M_H 
\end{align*}
By definition they are isomorpisms precisely for the $M$-reflexive modules 
\cite{A-F}. Either one of the $M$-dual functors map relexive modules to 
reflexive ones, so the restriction of $J$ and $K$ to the $M$-reflexive modules 
provides an adjoint equivalence 
\begin{equation} \label{eq: adj equiv}
\M_H^{\Mref}\sim(\,_N\M_N^\Mref)^\op\,,
\end{equation}
that is to say, a duality between the reflexive modules themselves.

Since $M$ has monoid structures both in $\M_H$ and $_N\M_N$, the functor $J$ is 
monoidal and $K$ is opmonoidal, 
\begin{align*}
\left.
\begin{alignedat}{2}
J_{X,X'}&:JX\oR JX'\to J(X\oN X')&\qquad (\xi\oR \xi')&\mapsto 
\hat\mu\circ(\xi\oN\xi')\\
J_0&:R\to JN&\qquad r&\mapsto\{n\mapsto \hat\eta(n)\ract r\}\\
\end{alignedat}
\right\}&\in \M_H\\
\left.
\begin{alignedat}{2}
K^{Y,Y'}&:KY\oN KY'\to K(Y\oR Y')&\qquad (\beta\oN\beta')&\mapsto
\mu\circ(\beta\oR\beta')\\
K^0&:N\to KR&\qquad n&\mapsto\{r\mapsto n\cdot\eta(r)\}
\end{alignedat}
\right\}&\in \,_N\M_N
\end{align*}
They are mates under the given adjunction $K\dashv J$, that is to say,
\begin{align}
K^{JX,JX'}\circ(\sigma_X\oN\sigma_{X'})&=KJ_{X,X'}\circ\sigma_{X\oN X'}\\
K^0&=KJ_0\circ\sigma_N\\
J_{KY,KY'}\circ(\sigma_Y\oR\sigma_{Y'})&=JK^{Y,Y'}\circ\sigma_{Y\oR Y'}\\
J_0&=JK^0\circ\sigma_R
\end{align}
These equations are simple consequences of the fact that the two monoid 
structures on $M$ come from the same $k$-algebra structure,
\[
\begin{CD}
M\o M@>>>M\oR M\\
@VVV @VV{\mu}V\\
M\oN M@>\hat\mu>>M
\end{CD}
\qquad 
\begin{CD}
k@>>> R\\
@VVV @VV{\eta}V\\
N@>\hat\eta>>M
\end{CD}
\]
We have, as in Theorem \ref{thm: U} (2), that $K$ is normal iff the map $N\to 
 M^H $ is an isomorphism and $J$ is normal iff the map $R\to M^N$ is 
an isomorphism.

In order to find monoidal subcategories in $\M_H$ and $_N\M_N$ that become 
monoidally dual under (\ref{eq: adj equiv}) we have to make further assumptions. 
Assume that the right bialgebroid $H$ is that of the horizontal Hopf algebroid 
of a distributive double algebra $A$ and assume that $N\subset M$ is $A$-Galois. 
We know from Lemma \ref{lem: fgp mon} and Theorem \ref{thm: U} that $\M_H^\fgp$ 
is a full  monoidal subcategory and the restriction of $K$ to this subcategory 
is strong opmonoidal. Therefore the restriction of $K$ will provide a monoidal 
equivalence iff all the fgp $H$-modules are $M$-reflexive. Since the class of 
reflexive modules is closed under taking direct summands and finite direct sums, 
this happens precisely when the regular object $H_H$ is $M$-reflexive. 
\[
\sigma_H:H\to\Hom_{N^e}(\Hom_H(H,M),M)\iso\End\,_NM_N
\]
being just the canonical embedding to the endomorphism Hopf algebroid $E$ we are 
left with considering the case when $H$ is $E$ and acts canonically on $M$.

\begin{thm} \label{thm: qsym}
Let $E$ be the endomorphism Hopf algebroid associated to the balanced 
depth 2 Frobenius extension $N\subset M$. Then the functor
\[
\Hom_E(\under, M):\M_E^\fgp\to \,_N\M_N^\op
\]
provides a monoidal duality between the categories of all fgp $E$-modules and 
those $N$-$N$-bimodules that are direct summands of  finite direct sums of 
$M$'s.
\end{thm}
\begin{proof}
$E_E$ is $M$-reflexive by construction. Thus $\M_E^\fgp$ is 
a full subcategory of the category of reflexive modules and (\ref{eq: adj 
equiv}) restricts to a category equivalence $F$. Since $\M_E^\fgp$ is generated 
by direct sums and direct summands from $E_E$, the same holds for the 
$N$-$N$-bimodules in the image of $F$ and for $F(E_E)\cong M$. The 
extension $N\subset M$ is $E$-Galois therefore $\Hom_E(\under,M)$ is strong 
monoidal on fgp modules. 
\end{proof}

Depending on the applications the content of the theorem varies from 
trivialities to nontrivial statements. For example, if $k\subset K$ is a 
separable field extension - including the case of classical Galois field 
extensions - then $E=\End K_k$ is the Hopf algebroid version of the weak Hopf 
algebra constructed in \cite{Sz: Strasbourg} and its representation category is 
trivial: The theorem reduces to the statement that the category of finite 
dimensional $k$-vector spaces is selfdual. 

If $N$ is a strongly $G$-graded $k$-algebra for a finite group $G$ and 
$N_e\subset N$ has centralizer $k\cdot 1$ then choosing 
$M=N\mash (kG)^*$ we obtain that $E$ is the group algebra $kG$ acting
canonically on the smash product $M$.

If $N$ is the observable algebra in rational quantum field theory and $M$ is 
the algebra of charge creating fields then $_NM$ is freely generated by 
finitely many fields $f_q^i\in M$ each of them implementing a localized 
endomorphism $\rho_q$ of $N$, i.e., $f^i_qn=\rho_q(n)f_q^i$, $n\in N$, 
$i=1,\dots, I_q$. The Doplicher-Haag-Roberts category $\DHR(N)$ is the full 
subcategory of $\End N$ the objects of which are finite direct sums of 
$\rho_q$'s and is a monoidal category by composition of endomorphisms.
One has a contravariant monoidal equivalence between $\DHR(N)$ and the category 
of $N$-$N$-bimodules that are finite direct sums of the bimodules $Nf_q^i$. 
Hence Theorem \ref{thm: qsym} gives a monoidal equivalence $\M_E^\fgp\simeq 
\DHR(N)$ and therefore the Hopf algebroid $E$ can be interpreted as the 
global gauge symmetry of the superselection sectors.

\end{document}